\documentclass[draft]{amsart}

\usepackage{graphics}
\usepackage[all]{xy}

\usepackage{amssymb,amsmath}
\usepackage{psfrag}

\usepackage{stmaryrd}

\usepackage{amsfonts}
\usepackage{amscd}

\usepackage{graphicx}
\usepackage{epstopdf}
\usepackage{tikz}
\usetikzlibrary{arrows,automata}

\usepackage{graphicx,epsfig,float}

\usepackage{url}

\makeatletter
\def\url@leostyle{%
  \@ifundefined{selectfont}{\def\UrlFont{\sf}}{\def\UrlFont{\small\ttfamily}}}
\makeatother
\newtheorem{thm}{Theorem}[section]
\newtheorem{lem}[thm]{Lemma}
\newtheorem{prop}[thm]{Proposition}
\newtheorem{cor}[thm]{Corollary}
\newtheorem{fact}[thm]{Fact}
\newtheorem{clm}[thm]{Claim}

\newtheorem{defn}[thm]{Definition}
\newtheorem{nrmk}[thm]{Remark}
\newtheorem{expl}[thm]{Example}

\newcommand{\pf}{{\bf Proof. }}

\renewcommand{\tilde}{\widetilde}
\renewcommand{\bar}{\overline}

\newcommand{\bB}{{\mathbf B}}

\renewcommand{\mod}{\mathrm{Mod}}

\newcommand{\op}{\mathrm{Op}}
\newcommand{\opc}{\mathrm{Op}^{{\rm cons}}}

\newcommand{\df}{\mathrm{def}}
\newcommand{\Df}{\mathrm{Def}}

\newcommand{\tDf}{\tilde{\mathrm{Def}}}
\newcommand{\tDfs}{\tilde{\mathrm{Def}(\mathbb S)}}
\newcommand{\tDfw}{\bigwedge \textrm{}\tilde{\mathrm{Def}}}

\newcommand{\pt}{\mathrm{pt}}

\begin{document}

\title {On definably proper maps}

\author {M\'{a}rio J. Edmundo}

\address{ Universidade Aberta\\ Rua Braamcamp 90\\
   1250-052 Lisboa, Portugal\\
 and \\
 CMAF Universidade de Lisboa\\
Av. Prof. Gama Pinto 2\\
1649-003 Lisboa, Portugal}

\email{mjedmundo@fc.ul.pt}

\author{Marcello Mamino}

\address{Laboratoire d'Informatique de l'\'Ecole Polytechnique (LIX)\\
B\^atiment Turing, bureau 2011\\
1 Rue Honor\'e d'Estienne d'Orves, Campus de l'\'Ecole Polytechnique\\
91120 Palaiseau, France}

\email{mamino@lix.polytechnique.fr}

\author{Luca Prelli}

\address{ CMAF Universidade de Lisboa\\
Av. Prof. Gama Pinto 2\\
1649-003 Lisboa, Portugal}


\email{lmprelli@fc.ul.pt}

\date{\today}
\thanks{The first author was supported by Funda\c{c}\~ao para a Ci\^encia e a Tecnologia, Financiamento Base 2008 - ISFL/1/209.  The third author is a member of the Gruppo Nazionale per l'Analisi Matematica, la Probabilit\`a e le loro Applicazioni (GNAMPA) of the Istituto Nazionale di Alta Matematica (INdAM) and was supported by Marie Curie grant PIEF-GA-2010-272021. This work is part of the FCT projects PTDC/MAT/101740/2008 and PTDC/MAT/122844/2010.\newline
 }


\begin{abstract}
In this paper we work in o-minimal structures with definable Skolem functions and show that: (i) a Hausdorff definably compact definable space is definably normal; (ii)  a continuous definable map between  Hausdorff locally definably compact definable spaces is definably proper if and only if it is  a proper morphism in the category of definable spaces. We  give several other characterizations of definably proper including one involving the existence of limits of definable types. We  also prove the basic properties of definably proper maps and  the invariance of definably proper (and definably compact) in elementary extensions and o-minimal expansions. 
\end{abstract}

\maketitle

\begin{section}{Introduction}\label{section intro}
Let ${\mathbb M}=(M,<, \ldots )$ be an arbitrary  o-minimal structure with definable Skolem functions. In this paper we  show that  Hausdorff definably compact definable spaces are definably normal (Theorem \ref{thm def comp hausd normal} ). 
We also show a local almost everywhere curve selection for Hausdorff locally definably compact definable spaces (Theorem \ref{thm local aecs}). 

Theorem \ref{thm def comp hausd normal} was only known in special cases: it was proved by Berarducci and Otero for definable manifolds in o-minimal expansions of real closed fields (\cite[Lemma 10.4]{bo01} - the proof there works as well in o-minimal expansions of ordered groups); it was proved in \cite{et} for definably compact groups in arbitrary o-minimal structures. Theorem \ref{thm local aecs} is an extension of the  almost everywhere curve selection for $M^n$ in arbitrary o-minimal structures proved by Peterzil and Steinhorn (\cite[Theorem 2.3]{ps}).

In Corollary \ref{cor def proper inv1} and Proposition \ref{prop def comp and s comp}  we show that definably compact is invariant under elementary extensions and  o-minimal expansions of ${\mathbb M}.$ In Proposition  \ref{prop def comp and comp} we show that if ${\mathbb M}$ is an o-minimal expansion of the ordered set of real numbers, then definably compact  corresponds to  compact. These invariance and comparison results extend similar results for definably compact subsets of $M^n$ in arbitrary o-minimal structures and answer partially a question from \cite{ps}.\\


In the authors recent work on the  formalism of the six Grothendieck operations on o-minimal sheaves (\cite{ep3} and \cite{ep4}) we require the basic theory  of morphisms  proper  in the category of o-minimal spectral spaces similar to the theory of proper morphisms in semi-algebraic geometry (\cite[Section 9]{dk1}) (and also  in algebraic geometry \cite[Chapter II, Section 4]{Har} or \cite[Chapter II, Section 5.4]{ega2}). Here, in Section \ref{section  def proper maps}, we provide such a theory by giving a category theory characterization of definably proper maps (as separated and universally closed morphisms in the category of definable spaces) and by proving the basic properties of such morphisms.   

In Theorems \ref{thm def proper inv1} and \ref{thm def proper inv3} we show that definably proper is invariant under elementary extensions and o-minimal expansions of ${\mathbb M}.$ In Theorem \ref{thm def proper inv4} we show that if ${\mathbb M}$ is an o-minimal expansion of the ordered set of real numbers, then definably proper corresponds to proper. These invariance and comparison results transfer to the notion of proper morphism in  the category of o-minimal spectral spaces.

The  formalism of the six Grothendieck operations on o-minimal sheaves (\cite{ep3} and \cite{ep4}) provides the cohomological ingredients required for  the computation of the subgroup of $m$-torsion points of a definably compact, abelian definable group $G$  - extending the main result of  \cite{eo} which was proved in o-minimal expansions of ordered fields using the o-minimal singular (co)homology. This result is enough to settle  Pillay's conjecture for definably compact definable groups (\cite{p2} and \cite{hpp1})
in arbitrary o-minimal structures. See \cite{oh5p}. Pillay's conjecture is a non-standard analogue of Hilbert's 5$^{{\rm o}}$ problem for locally compact topological groups, roughly it says that after taking the quotient by a ``small subgroup'' (a smallest type-definable subgroup of bounded index)  the quotient when equipped with the the so called logic topology is a compact real Lie group of the same dimension.\\

Finally in Section \ref{section def comp, def proper and def types}  we  prove that definable compactness of Hausdorff definable spaces can be characterized by the existence of limits of definable types (Theorem \ref{thm def comp norm types}), extending a remark by Hrushovski and Loeser (\cite{HrLo}) in the affine case. In Theorem \ref{thm def proper1} we prove a corresponding characterization of definably proper maps between Hausdorff locally definably compact definable spaces which, when transferred to  morphisms proper  in the category of o-minimal spectral spaces, is the analogue  of the valuative criterion  for properness in algebraic geometry (\cite[Chapter II, Theorem 4.7]{Har}). As it is known, in o-minimal structures with definable Skolem functions, definable types correspond to valuations (\cite{MarStein} and \cite{p94}).

\medskip
\emph{Acknowledgements.} We wish to thank the referee for reading our paper so carefully  and for making  very detailed reports with suggestions of simplifications of some proofs as well as fixes to gaps in previous versions of the paper. \\

\end{section}

\begin{section}{On definably compact  spaces}\label{section def compact haus}

\begin{subsection}{Hausdorff definably compact spaces}\label{subsection Hausd def comp}
Here we will show that if ${\mathbb M}$ has definable Skolem functions, then Hausdorff definably compact definable spaces are definably normal.\\

Below we will assume the reader familiarity with basic o-minimality (see for example \cite{vdd}). Below,  by definable we mean definable in ${\mathbb M}$ possibly with parameters. Recall also that ${\mathbb M}$ has definable Skolem functions if for every uniformly definable family $\{F_t\}_{t\in T}$ of definable sets, then there is a definable map $f:T\to \bigcup _tF_t$ such that for each $t\in T$ we have $f(t)\in F_t.$\\

Recall the notion of definable spaces (\cite{vdd}):\\

\begin{defn}\label{defn def space}
{\em
A {\it definable space} is a tuple $(X,(X_i,\theta _i)_{i\leq k})$  where:

 \begin{itemize}
\item[$\bullet  $]
$X=\bigcup _{i\leq k}X_i$;

\item[$\bullet  $]
each $\theta  _i:X_i\rightarrow  M^{n_i}$ is an injection such that $\theta _i(X_i)$ is a definable subset of $M^{n_i}$ with the induced topology;

\item[$\bullet  $]
for all $i, j$, $\theta  _i(X_i\cap X_j)$ is an open definable subset of  $\theta _i(X_i)$ and  the transition maps $\theta _{ij}:\theta  _i(X_i\cap X_j)\rightarrow  \theta  _j(X_i\cap X_j):x\mapsto \theta  _j(\theta  _i^{-1}(x))$ are definable homeomorphisms.
\end{itemize}
We call the $(X_i, \theta _i)$'s the {\it definable charts of $X$} and define the dimension of $X$ by  $\dim X=\max \{\dim \theta_i(X_i):i=1,\dots ,k\}$. If all the $\theta _i(X_i)$'s are open definable subsets of some $M^n$, we say that $X$ is a {\it definable manifold of dimension $n$}.

A  definable space $X$ has a topology such that each $X_i$ is open and the $\theta  _i$'s are homeomorphisms: a subset $U$ of $X$ is  open in the basis of this topology if and only if for each $i$, $\theta  _i(U\cap X_i)$ is an open definable subset of $\theta  _i(X_i)$.

A map $f:X\to Y$ between  definable spaces with definable charts $(X_i,\theta _i)_{i\leq k}$ and $(Y_j,\delta _j)_{j\leq l}$ respectively is {\it a  definable map} if:
\begin{itemize}
\item[$\bullet$]  for each $i$ and every $j$ with $f(X_i)\cap Y_j\neq\emptyset $, $\delta _j\circ f\circ \theta _i^{-1}:\theta _i(X_i)\to \delta _j(Y_j)$ is a definable map between definable sets.
\end{itemize}
 We say that a definable space is {\it affine} if it is definably homeomorphic to a definable set with the induced topology.
 }
 \end{defn}

The construction above defines  the {\it category of  definable spaces with  definable continuous maps} which we denote by $\Df$. All topological notions on definable spaces are relative to the topology above. Note however, that often we will have to replace topological notions on definable spaces by their definable analogue.

We  say that a subset $A$ of a definable space $X$ is definable if and only if for each $i$,  $\theta _i(A\cap X_i)$ is a definable subset of $\theta  _i(X_i)$.   A definable subset $A$ of a definable space $X$ is naturally a definable space and its topology is the induced topology, thus we also call them {\it definable subspaces}.  \\

 In nonstandard o-minimal structures closed and bounded definable sets are not compact. Thus we have to replace the notion of compactness by a suitable definable analogue.\\

Let $X$ be a definable space and $C\subseteq X$ a definable subset. By a {\it definable curve in $C$} we mean a continuous definable map  $\alpha :(a,b)\to C\subseteq X$, where $a<b$ are in $M\cup \{-\infty, +\infty \}.$ We say that a definable curve $\alpha :(a,b)\to C\subseteq X$ in $C$ is {\it completable in C} if both limits  $\lim _{t\to a^+}\alpha (t)$ and $\lim _{t\to b^-}\alpha (t)$ exist in $C,$  equivalently if there exists a continuous definable map  $\bar{\alpha }:[a,b]\to C\subseteq X$ such that the diagram
$$
\xymatrix{
(a,b)  \ar[r]^{\alpha } \ar@{^{(}->}[d] & C\subseteq X\\
[a,b]  \ar@{-->}[ur]_-{\overline{\alpha }} & }
$$
is commutative. \\

\begin{defn}\label{defn def compact}
{\em
Let $X$ be a definable space and $C\subseteq X$ a definable subset. We say that $C$ is {\it definably compact} if every definable curve in $C$ is completable in $C$ (see \cite{ps}).\\
}
\end{defn}

The following is easy:\\

\begin{fact}\label{fact def skolem def comp1}
Suppose  that ${\mathbb M}$ has  definable Skolem functions.
Let   $f:X\to Y$ a continuous definable map between definable spaces. If $K\subseteq X$ is a definably compact definable subset, 
then $f(K)$ is a definably compact definable subset of $Y$. \\
\end{fact}

For  definable subsets $X\subseteq M^n$ with their induced topology (i.e.  affine definable spaces) the notion of definably compact is very well behaved. Indeed, we have (\cite[Theorem 2.1]{ps}):

\begin{fact}\label{fact def comp affine}
A definable subset $X\subseteq M^n$ is definably compact if and only if it is closed and bounded in $M^n$  \\ 
\end{fact}

However, in general,  
definably compact definable subsets of a definable space are not Hausdorff and are not even necessarily  closed subsets:

\begin{expl}[Non Hausdorff and non closed definably compact subsets]\label{expl def comp no haus, closed}
{\em







Let $a, b, c, d\in M$ be such that  $c<a<b<d.$ Let $X$ be the definable space with definable charts $(X_i, \theta _i)_{i=1,2}$ given by:
 $X_1= (\{\langle x, y \rangle \in [c,d]\times [c, d]: x=y\}\setminus  \{\langle b, b\rangle \})\cup \{\langle b, a\rangle \}\subseteq  M^2$, $X_2=\{\langle x, y \rangle \in [c,d]\times [c, d]: x=y\}\subseteq M^2$ and  $\theta _i=\pi _{|X_i}$ where $\pi :M^2\to M$ is the projection onto the first coordinate. Then any open definable neighborhood in $X$ of the point $\langle b, a\rangle$  intersects   any open definable neighborhood in $X$ of the point $\langle b, b\rangle.$ Clearly $X$ is definably compact but not Hausdorff and $X_2$ is a definably compact subset which is not closed (in $X$).

 }
 \end{expl}

It is desirable to work in a situation where definably compact subsets are closed. We will show that this is the case in Hausdorff definable spaces when ${\mathbb M}$ has definable Skolem functions. \\

Before we need to  introduce some notations.\\

 Let $X$ be a definable space and let $(X_i, \theta _i)_{i\leq k}$ be the definable charts of $X$ with $\theta _i(X_i)\subseteq M^{n_i}.$ Let $N=n_1+\cdots + n_k$ and fix  a point $*\in M.$ For each $i\leq k,$ let $\pi _i:M^N=M^{n_1}\times \cdots \times  M^{n_k}\to M^{n_i}$ be the natural projection and let $\rho _i:M^{n_i}\to M^{N}$ be the  inclusion with $\rho _i(M^{n_i})=\{*\}\times \cdots \times \underbrace{M^{n_i}}_\text{position $i$}\times \cdots \times \{*\}\subseteq M^N.$ Identify each $M^{n_i}$ with $\rho _i(M^{n_i})\subseteq M^N.$  Identify as well each $\theta _i(X_i)$ with  $\rho _i(\theta _i(X_i))\subset M^N$ and each $\theta _i$ with $\rho _i\circ \theta _i.$\\

 For $a\in X$ let $I_a=\{i\leq k: a\in X_i\}$ and set 
{\small
 $$D(a)=\{\langle \langle d^-_1,d^+_1\rangle , \ldots ,\langle d^-_{N}, d^+_{N}\rangle \rangle \in M^{2N}: \theta _j (a)\in \pi _j(\Pi _{l=1}^{N}(d^-_l,d^+_l)) \,\,\textrm{for all}\,\, j\in I_a\}.$$  
 }

 Consider the finite set $I_X=\{I\subseteq \{1,\ldots , k\}: I=I_a \,\,\textrm{for some}\,\, a\in X\}.$ Then   each $X_I=\{x\in X: I_x=I\}$ with $I\in I_X$ is a definable subset and $X=\bigsqcup _{I\in I_X}X_I.$  
 Therefore, 
 $$\{ D(a)\}_{a\in X}$$
  is a uniformly definable family of definable sets, since it  is defined by the first-order formula
 $$\bigvee _{I\in I_X}[(a\in X_I)\land \bigwedge _{j\in I}\bigwedge _{l=N_{j-1}+1}^{N_j}(d^-_{l}<\theta _j(a)_l<d^+_{l})]$$
 where for each $i\leq k$  we set $N_i=n_1+\cdots +n_{i}$ and where $\theta _j(a)_l$ is the $l$-coordinate of $\theta _j(a)\in M^N.$

 For $d, d' \in D(a)$ we set $d\preceq d'$ if   we have $\Pi _{l=1}^{N}(d^-_{l},d^+_{l})\subseteq \Pi _{l=1}^{N}(d'^-_{l},d'^+_{l})$ and the notation  $d\prec d'$ is used whenever   we have $\Pi _{l=1}^{N}(d^-_{l},d^+_{l})\subset \Pi _{l=1}^{N}(d'^-_{l},d'^+_{l}).$ \\
 
 The following are immediate:\\
 
 \begin{itemize}
 \item[(D0)] 
 The relation $\preceq $ on $D(a)$ is a definable downwards directed order  on $D(a).$
 \item[(D1)] 
 The set $D(a)\subseteq M^{2N}$ is an open definable subset.
 \item[(D2)]
 If $d\in D(a)$ then $\{d'\in D(a):d'\prec d\}$ is an open definable subset of $D(a).$\\
\end{itemize}
 
 For $a\in X$ and $d=\langle \langle d^-_1,d^+_1), \ldots , \langle d^-_{N}, d^+_{N}\rangle \rangle \in D (a)$  set 
$$U(a, d)=\bigcap _{j\in I_a}\theta _j ^{-1}(\theta _j(X_j)\cap \pi _j(\Pi _{l=1}^{N}(d^-_{l},d^+_{l}))).$$
 Then  
 $$\{U(a,d)\}_{d\in D(a)}$$
  is a uniformly definable system of fundamental open definable neighborhoods of $a$ in $X.$ \\

The following will also be useful:\\

\begin{itemize}
\item[(D3)]
If $a, a'\in X$ are such  that $I_{a'}\subseteq I_a$, then for every $d\in D(a)\cap D(a')$ we have $U(a,d)\subseteq U(a',d).$\\
\end{itemize}

Finally we will also require:\\

\begin{itemize}
\item[(D4)]
If $a\in X$  and $W$ is an open definable neighborhood of $a$ then the set $\{d\in D(a): U(a,d)\subset W\}$ is an open definable subset of $D(a).$\\
\end{itemize}

If $B\subseteq X$ is a definable subset and $\epsilon :B\to M^{2N}$ is a definable map such that $\epsilon (x)\in D(x)$ for all $x\in B$, then
$$U(B, \epsilon )=\bigcup _{x\in B}U(x,\epsilon (x))$$
is an open definable neighborhood of $B$ in $X.$\\

It follows that: 

\begin{nrmk}\label{nrmk closed open of}
{\em
The notions of open (resp. closed) in a definable space $X$ are first-order in the sense that if $(A_t)_{t\in T}$ is a uniformly definable family of definable subsets of $X$, then the set of all $t\in T$ such that $A_t$ is an open (resp. a closed) subset of $X$ is a definable set.\\
}
\end{nrmk}

Recall that a topological space $X$ is \textit{regular} if one the following  equivalent conditions holds:
\begin{enumerate}
\item
for every $a\in X$ and $S\subseteq X$ closed such that $a\not \in S$, there are open disjoint subsets $U$ and $V$ of $X$ such that $a\in U$ and $S\subseteq V$;
\item
for every $a\in X$ and $W\subseteq X$ open such that $a\in W$, there is $V$ open subset of $X$ such that $a\in V$ and $\bar{V}\subseteq W$.\\
\end{enumerate}

\begin{prop}\label{prop  def haus reg}
Suppose  that ${\mathbb M}$ has  definable Skolem functions. Let $X$ be a Hausdorff  definable space. Then for any $a\in X$ and any definably compact subset $K\subseteq X$ such that $a\not \in K,$ there are finitely many  definably compact subsets   $K_i$ ($i=1, \ldots , l$)  of $K,$ finitely many continuous definable functions $\epsilon _i:K_i \to M^{2N}$  with $\epsilon _i(x)\in D(x)$ for all $x\in K_i$  and there is $d\in D(a)$ such that:
\begin{itemize}
\item
$K\subseteq \bigcup _{i=1}^lU(K_i, \epsilon _i).$
\item
 $U(a,d)\cap (\bigcup _{i=1}^lU(K_i, \epsilon _i))=\emptyset .$
\end{itemize} 
In particular, if $X$ is a Hausdorff, definably compact  definable space, then $X$  is regular. \\
\end{prop}

\pf
We fix $a\in X$ and prove the result by induction on the dimension of definably compact subsets $K\subseteq X$ such that $a\not \in K.$ 

If $\dim K=0$, then this follows because $X$ is Hausdorff. Assume the result holds for every definably compact subset $L$ of $X$ such that $a\not \in  L $ and $\dim L<\dim K.$

Since $X$ is Hausdorff, for each  $x\in K$ there is $d'\in D(x)$ and there is $d\in D(a)$ such that $U(a,d)\cap U(x, d')=\emptyset .$ By definable Skolem functions there  are definable maps
$$g: K \to M^{2N}$$
and 
$$h: K\to M^{2N}$$
such that: 
\begin{itemize}
\item[(i)]  
$g(x)\in D(a)$ for all $x\in K;$ 
\item[(ii)]    $h(x)\in D(x)$ for all $x\in K;$ 
\item[(iii)]
$U(a,g(x))\cap U(x, h(x))=\emptyset .$ 
\end{itemize}
Since, by Remark \ref{nrmk closed open of}, continuity is first-order,  the  subset of $K$ where either $g$ or $h$ is not continuous is a definable subset. By working in charts and using \cite[Chapter 3, (2.11) and Chapter 4, (1.8)]{vdd} this definable subset has dimension $<\dim K$ and, if  $L$ is the closure of this subset, then $\dim L<\dim K.$  By induction hypothesis,  there are finitely many  definably compact subsets   $L_i$ ($i=1, \ldots , k$)  of $L,$ finitely many continuous definable functions $\epsilon _i:L_i \to M^{2N}$  with $\epsilon _i(x)\in D(x)$ for all $x\in L_i$  and there is $d_L\in D(a)$ such that:
\begin{itemize}
\item
$L\subseteq \bigcup _{i=1}^kU(L_i, \epsilon _i).$
\item
 $U(a,d_L)\cap (\bigcup _{i=1}^kU(L_i, \epsilon _i))=\emptyset .$
\end{itemize}

Let $K'=K\setminus \bigcup _{i=1}^kU(L_i,\epsilon _i).$ Then $K'$ is definably compact and both $g_{|}: K' \to M^{2N}$ and  $h_{|}:K'\to M^{2N}$ are continuous.  We show that there is $d_{K'}\in D(a)$ such that $d_{K'}\preceq g_{|}(x)$ for all $x\in K'.$ 

Write $g(x)=\langle \langle g^-(x)_{1},g^+(x)_{1}\rangle, $ $ \ldots ,$ $\langle g^-(x)_{N}, $ $ g^+(x)_{N} \rangle \rangle \in D(a)$  where for  each $l=1, \ldots , N$, $g^-(x)_l$ and $g^+(x)_l$ are the two $l$-components of $g(x).$ 
By Fact \ref{fact def skolem def comp1}, for each $l=1,\ldots , N,$ let $d^-_{l}=\max\{g^-(x)_l:x\in K'\}$ and $d^+_{l}=\min\{g^+(x)_l:x\in K'\}.$ Since each $d^-_l=g^-(z)_l$ for some $z\in K',$ and similarly  each $d^+_l=g^+(z')_l$) for some $z'\in K',$  we have $d_{K'}:=\langle \langle d^-_{1},d^+_{1}\rangle , \ldots , \langle d^-_{N}, d^+_{N}\rangle \rangle \in D(a).$ By construction  we also have $d_{K'}\preceq g_l(x)$ for all $x\in K'.$  


To finish the proof,  choose $d\preceq d_L, d_{K'}$ by (D0) and, for each $i=1,\ldots , k,$ set $K_i=L_i$  and take also $K_{k+1}=K'$ and $\epsilon _{k+1}=h_{|K'}.$ Then, by construction, 
\begin{itemize}
\item
$K\subseteq \bigcup _{i=1}^{k+1}U(K_i, \epsilon _i).$
\item
 $U(a,d)\cap (\bigcup _{i=1}^{k+1}U(K_i, \epsilon _i))=\emptyset .$
\end{itemize} 
\qed \\

 The following is now immediate:

\begin{cor}\label{cor def comp2}
Suppose  that ${\mathbb M}$ has  definable Skolem functions. Suppose that $X$ is a Hausdorff definable space. If $K$ is definably compact subset of $X$, then $K$ is a closed definable subset.\\
\end{cor}


We will require the following:

\begin{lem}\label{lem closure}
Suppose  that ${\mathbb M}$ has  definable Skolem functions. Let $X$ be a Hausdorff, definably connected,  definable space and $K\subseteq X$  a definably compact subset. Let $\epsilon :K\to M^{2N}$ be a  definable continuous map such that $\epsilon (x)\in D(x)$ for all $x\in K$ and suppose that for each $w\in K$ there is $d\in D(w) $ such that $\epsilon (w)\prec d$ and  $\bar{U(w, d)}$ is definably compact. Then 
$$\bigcup _{x\in K}\bar{U(x, \epsilon (x))}$$
is a closed definably compact definable neighborhood of $K$. In particular we have
$$\bar{U(K,\epsilon)}=\bar{\bigcup _{x\in K}U(x, \epsilon (x))}=\bigcup _{x\in K}\bar{U(x, \epsilon (x))}.$$
\end{lem}

\pf
Let $\alpha : (a,b)\to \bigcup _{x\in K}\bar{U(x, \epsilon (x))}$ be a definable curve. We have show that the limit $\lim_{t\to b^-}\alpha (t)$ exists in $\bigcup _{x\in K}\bar{U(x, \epsilon (x))}.$

By definable Skolem functions there is a definable map  $\beta :(a,b)\to K$ such that for each $t\in (a,b)$ we have 
$$\alpha (t)\in \bar{U(\beta (t), \epsilon (\beta (t)))}.$$ 
By o-minimality, after shrinking $(a,b)$ if necessary,  i.e., after replacing $a$ by $a'\in (a,b)$ if needed, we may assume that $\beta $ is a definable curve in $K.$ Since $K$ is definably compact, let $w=\lim _{t\to b^-}\beta (t)\in K.$ 
 Let also $\bar{\beta }:(a,b]\to K$ be the continuous definable map such that $\bar{\beta }_{|(a,b)}=\beta _{|(a, b)}.$  


Recall that we have $\epsilon \circ \bar{\beta } (b)=\epsilon (w)\in D(w)$ and $D(w)\subseteq M^{2N}$ is an open definable subset by (D1). 
Since $\epsilon :K\to M^{2N}$ is continuous, it follows from the continuity of $\epsilon \circ \bar{\beta }: (a,b]\to  M^{2N}$ at $b$ that there is $a'\in (a, b)$ 
such that $\epsilon \circ \bar{\beta }(t)\in D(w)$ for all $t\in [a', b]$. 

Since for each $j\in I_w$, $X_j$ is an open definable neighborhood of $w$, by continuity, after shrinking $(a',b]$ if necessary, we may assume that $\bar{\beta }(t)\in X_j$ for all $t\in [a',b]$ and all $j\in I_w.$ 
Thus we must have $I_w\subseteq I_{\bar{\beta }(t)}$ for all  $t\in [a',b].$ Therefore,  by (D3), for all $t\in [a', b]$ we have $U(\bar{\beta }(t),\epsilon (\bar{\beta }(t)))\subseteq U(w, \epsilon (\bar{\beta }(t))).$ 

In particular, for each $t\in [a',b)$ we have 
$$\alpha (t)\in \bar{U(w, \epsilon (\bar{\beta }(t)))}.$$ 

By hypothesis there is $d\in D(w) $ such that $\epsilon (w)=\epsilon (\bar{\beta  }(b))\prec d$ and  $\bar{U(w, d)}$ is definably compact. By (D2) and continuity of $\epsilon \circ \bar{\beta }:[a',b]\to D(w)\subseteq M^{2N}$, after shrinking $(a',b]$ if necessary, we may further  assume that $\epsilon (\bar{\beta  }(t))\prec d$ for all $t\in [a', b].$ Therefore, $$\alpha (t)\in \bar{U(w, d)}$$ for all $t\in [a', b).$ 

Since $\bar{U(w, d)}$ is definably compact, there exists the limit $\lim _{t\to b^-}\alpha (t) \in \bar{U(w, d)}.$ Let $v=\lim _{t\to b^-}\alpha (t) \in \bar{U(w, d)}.$ We want to show that   $v\in \bar{U(w, \epsilon (w))}.$ Suppose not and set $L=\bar{U(w, \epsilon (w))}.$ Since $L$ is definably compact subset of $\bar{U(w, d)}$, by Proposition \ref{prop  def haus reg}, there are finitely many  definably compact subsets   $L_i$ ($i=1, \ldots , k$)  of $L,$ finitely many continuous definable functions $\epsilon _i:L_i \to M^{2N}$  with $\epsilon _i(x)\in D(x)$ for all $x\in L_i$  and there is $d_L\in D(v)$ such that:
\begin{itemize}
\item
$L\subseteq \bigcup _{i=1}^kU(L_i, \epsilon _i).$
\item
 $U(v,d_L)\cap (\bigcup _{i=1}^kU(L_i, \epsilon _i))=\emptyset .$
 \end{itemize}
We have $U(w,\epsilon (w))\subseteq L \subseteq \bigcup _{i=1}^kU(L_i, \epsilon _i).$ If $U(w,\epsilon (w))=\bigcup _{i=1}^kU(L_i, \epsilon _i)$ then $U(w,\epsilon (w))=L=\bar{U(w, \epsilon (w))}$ and so $U(w,\epsilon (w))$ is a closed and open definable subset of $X$. Since $X$ is definably connected we would have $U(w, \epsilon (w))=X$ and so 
$v\in \bar{U(w, \epsilon (w))}$ which is a contradiction.
 
Since $U(w, \epsilon (w))\subset \bigcup _{i=1}^kU(L_i, \epsilon _i)$ and $\bigcup _{i=1}^kU(L_i, \epsilon _i)$ is an open definable neighborhood of $w$, by (D4) there is $a''\in [a',b]$ such that $U(w, \epsilon (\bar{\beta }(t))) \subset \bigcup _{i=1}^kU(L_i, \epsilon _i)$ for all $t\in [a'', b].$ Therefore, for each $t\in [a'',b]$ we have 
$$\alpha (t)\in \bar{\bigcup _{i=1}^kU(L_i, \epsilon _i)}.$$ 
This implies that $v\in \bar{\bigcup _{i=1}^kU(L_i, \epsilon _i)}$ which contradicts the fact that $U(v,d_L)\cap \bigcup _{i=1}^kU(L_i, \epsilon _i)=\emptyset .$\\

By Corollary \ref{cor def comp2}, $\bigcup _{x\in K}\bar{U(x, \epsilon (x))}$ is closed and hence 
$$\bar{U(K,\epsilon)}=\bar{\bigcup _{x\in K}U(x, \epsilon (x))}=\bigcup _{x\in K}\bar{U(x, \epsilon (x))}.$$
\qed \\

Recall that a  definable space $X$  is \textit{definably normal} if one of the following equivalent conditions holds:
\begin{enumerate}
\item
for every disjoint closed definable subsets $Z_1$ and $Z_2$ of $X$ there are disjoint open definable subsets $U_1$ and $U_2$ of $X$ such that $Z_i\subseteq U_i$ for $i=1,2.$

\item
for every  $S\subseteq X$ closed definable  and $W\subseteq X$ open definable such that $S\subseteq W$, there is an open definable subsets $U$  of $X$ such that $S\subseteq U$ and $\bar{U}\subseteq W$.\\
\end{enumerate}

In general regular does not imply definably normal:

\begin{expl}[Regular non definably normal definable space]\label{expl reg no def norm}
{\em
Assume that ${\mathbb M}=(M,<)$ is a dense linearly ordered set with no end points. Let $a, b, c, d\in M$ be such that $c<a<b<d$ and let $X=(c,d)\times (c,d)\setminus \{\langle a, b\rangle \}.$ Since $X$ is affine it is regular. Note also that the only open definable subsets of $X$ are the intersections with $X$ of definable subsets of $M^2$ which are  finite unions of  non empty finite intersections  $W_1\cap \cdots \cap W_k$ where each $W_i$ is either an  open box in $M^2$, $\{\langle x, y \rangle \in M^2: x<y\}$ or $\{\langle x, y \rangle \in M^2: y<x\}.$

Let $C=\{\langle x, y\rangle \in X: x=a\}$ and let $D=\{\langle x, y \rangle \in X: y=b\}.$ Then $C$ and $D$ are closed disjoint definable subsets of $X$. However, by the description of the open definable subset of $X$, there are no open disjoint definable subsets $U$ and $V$ of $X$ such that $C\subseteq U$ and $D\subseteq V.$\\
}
\end{expl}

\begin{thm}\label{thm def comp hausd normal}
Suppose  that ${\mathbb M}$ has  definable Skolem functions. If $X$ is a Hausdorff, definably compact  definable space, then $X$ is definably normal. In fact,  for every $K\subseteq X$  closed definable subset and for every $V\subseteq X$ open definable subset, if $K\subseteq V,$ then  there are finitely many  definably compact subsets   $K_i$ ($i=1, \ldots , l$)  of $K$ and finitely many continuous definable functions $\epsilon _i:K_i \to M^{2N}$  with $\epsilon _i(x)\in D(x)$ for all $x\in K_i$   such that:
\begin{itemize}
\item
$K\subseteq \bigcup _{i=1}^lU(K_i, \epsilon _i).$
\item
 $\bigcup _{i=1}^l\bar{U(K_i, \epsilon _i)}\subseteq V.$
\end{itemize}

\end{thm}

\pf
Clearly we may assume that $X$ is definably connected and we can fix  $V\subseteq X$ an open definable subset.  We prove the result by induction on the  dimension of  closed definable subsets $K\subseteq X$  such that $K\subset V.$ 

If $\dim K=0$ then the result follows since $X$ is regular (Proposition \ref{prop  def haus reg}). So assume that the result holds for every closed definable subset $L$ such that $L\subseteq V$ and $\dim L<\dim K.$

Since $X$ is regular (Proposition \ref{prop  def haus reg}),  for each  $x\in K$ there is $d\in D(x)$  such that  $\bar{U(x,d)}\subseteq V.$ Since the property ``$d\in D(x)$ and $\bar{U(x,d)}\subseteq V$'' is first-order (Remark \ref{nrmk closed open of}), by definable Skolem functions, there is a definable map
$$\delta  : K\to M^{2N}$$
such that, for all $x\in K$: 
\begin{itemize}
\item[(i)]   
 $\delta (x)\in D(x);$  
\item[(ii)] 
$\bar{U(x,\delta (x))}\subseteq V.$ 
\end{itemize}
By definable Skolem functions again and by (D2), there is a definable map
$$\epsilon  : K\to M^{2N}$$
such that, for all $x\in K$: 
\begin{itemize}
\item[(i)]   
 $\epsilon (x)\in D(x);$ 
\item[(ii)]
$\epsilon (x)\prec \delta (x);$
\item[(iii)] 
$\bar{U(x,\epsilon (x))}\subseteq V.$ 
\end{itemize}
Since, by Remark \ref{nrmk closed open of}, continuity is first-order,  the  subset of $K$ where  $\epsilon $  is not continuous is a definable subset. By working in charts and using \cite[Chapter 3, (2.11) and Chapter 4, (1.8)]{vdd} this definable subset has dimension $<\dim K$ and, if  $L$ is the closure of this subset, then $\dim L<\dim K.$  By induction hypothesis,  there are finitely many  definably compact subsets   $L_i$ ($i=1, \ldots , k$)  of $L$ and finitely many continuous definable functions $\epsilon _i:L_i \to M^{2N}$  with $\epsilon _i(x)\in D(x)$ for all $x\in L_i$   such that:
\begin{itemize}
\item
$L\subseteq \bigcup _{i=1}^kU(L_i, \epsilon _i).$
\item
 $\bigcup _{i=1}^k\bar{U(L_i, \epsilon _i)}\subseteq V.$
\end{itemize}

Let $K'=K\setminus \bigcup _{i=1}^kU(L_i, \epsilon _i).$ Then $K'$ is a closed definable subset  and  $\epsilon '=\epsilon _{|}: K' \to M^{2N}$ is continuous.  Furthermore, for each $w\in K'$ there is $d=\delta (w)\in D(w) $ such that $\epsilon '(w)\prec d$ and  $\bar{U(w, d)}$ is definably compact. Therefore,  by Lemma \ref{lem closure}, we have $\bar{U(K',\epsilon ')}\subseteq V.$ 

For each $i=1,\ldots , k,$ set $K_i=L_i$  and take also $K_{k+1}=K'$ and $\epsilon _{k+1}=\epsilon '.$ Then, by construction, 

\begin{itemize}
\item
$K\subseteq \bigcup _{i=1}^{k+1}U(K_i, \epsilon _i).$
\item
 $\bigcup _{i=1}^{k+1}\bar{U(K_i, \epsilon _i)}\subseteq V.$
\end{itemize} 

\qed \\

Definable normality gives the shrinking lemma (compare with \cite[Chapter 6, (3.6)]{vdd}):\\

\begin{cor}[The shrinking lemma]\label{cor shrinking lemma}
Suppose  that ${\mathbb M}$ has  definable Skolem functions. Let $X$ be a Hausdorff definably compact definable  space. If $\{U_i:i=1,\dots ,n\}$ is a covering of $X$ by open definable subsets, then there are definable open subsets $V_i$ and definable closed subsets $C_i$ of $X$ ($1\leq i\leq n$) with $V_i\subseteq C_i\subseteq U_i$ and $X=\cup \{V_i:i=1,\dots, n\}$. \\
\end{cor}

 \end{subsection}
 
 \begin{subsection}{Local almost everywhere curve selection}\label{subsection local aecs}
 To prove our results about definably proper maps later we will need a local version of an extension to definable spaces of  the almost everywhere curve selection (\cite[Theorem 2.3]{ps}):

\begin{fact}\label{fact arcs ps} 
If $C\subseteq M^n$ is a definable subset   which is not closed, then there is a definable set $E\subseteq \bar{C}\setminus C$ such that  $\dim E<\dim(\bar{C}\setminus C)$ and  for every $x\in  \bar{C}\setminus (C\cup E)$ there is a definable curve in $C$ which has $x$ as a limit point.\\
\end{fact}

We say that {\it the almost everywhere curve selection holds for a   definable space $X$} if  for every   definable subset  $C\subseteq X$  which is not closed, there is a definable set $E\subseteq \bar{C}\setminus C$ such that  $\dim E<\dim(\bar{C}\setminus C)$ and  for every $x\in  \bar{C}\setminus (C\cup E)$ there is a definable curve in $C$ which has $x$ as a limit point. \\

For general definable spaces, even affine ones, even if ${\mathbb M}$ has definable Skolem functions, almost everywhere curve selection does not hold: 

\begin{expl}\label{expl no aecs}
{\em
$\,$
\begin{enumerate}
\item
In ${\mathbb M}=({\mathbb Q}, <),$ for the definable set $D=\{\langle x,y\rangle \in {\mathbb Q}^2:0<y<x\}$ there is no definable curve in $D$ with limit $d=\langle 0,0\rangle .$ (This example is from \cite{ps}).
\item
Let $\Gamma =({\mathbb R},<, 0,  -, +, (q)_{q\in {\mathbb Q}}).$ Let $\Gamma _0=\{0\}\times \Gamma ,$ $\Gamma _1=\{1\}\times \Gamma $  and let  $\infty $ be a new symbol such that $\langle 0, x\rangle <\infty <\langle 1, y\rangle $ for all $x, y\in {\mathbb R}.$ Let $M=\Gamma _0\cup \{\infty \}\cup \Gamma _1$ be equipped with the natural induced total order from $<.$ Let ${\mathbb M}$ be the structure obtained by putting on $\Gamma _0$ and on $\Gamma _1$ the induced structure from $\Gamma .$  Then ${\mathbb M}$ has definable Skolem functions (since each copy of $\Gamma $ has definable Skolem functions by \cite[Chapter 6, (1.2)]{vdd}). However, for the definable set $D=\{\langle \langle 0,x\rangle ,\langle 0, y\rangle \rangle \in M^2:x,y>0\}$ there is no definable curve in $D$ with limit $d=\langle \langle 0,0\rangle ,\infty \rangle .$ Indeed, any definable curve in $D$ will be definable in $\Gamma $ and so its graph will be a piecewise linear  subset of $D$ (\cite[Chapter I, (7.8)]{vdd}). By piecewise-linearity  there are no definable bijections between bounded and unbounded intervals and so  no definable curve in $D$ will  have $d=\langle \langle 0, 0\rangle , \infty \rangle $ for a limit point. (This example is essentially the same as the $\Gamma _{\infty } $ from \cite[Section 4.1]{HrLo} - the only difference is that we added a new copy of $\Gamma ,$ the $\Gamma _1,$ so that our $M$ has no endpoints).
\end{enumerate}
In both cases, if $X=D\cup \{d\}$ then almost everywhere curve selection does not hold for the definable space $X$ since in $X$ we have $\bar{D}\setminus D=\{d\}.$\\
}
\end{expl}

The almost everywhere curve selection fails for $X\subseteq M^2$ in Example \ref{expl no aecs}  because $X$ there is not a locally closed subset of $M^2:$\\

\begin{lem}\label{lem aecs loc closed}
Suppose that $X$ is a  definable space  and that the almost everywhere curve selection holds for $X.$  Then the almost everywhere curve selection holds for every locally closed definable subset  of $X$.  
\end{lem}

\pf
Let $Z$ be a closed definable subset of $X$ and let $C\subseteq Z$ be a definable subset which is not closed in $Z$. Since  $Z$ is closed, $\bar{C}={\rm cl}_Z(C)\subseteq  Z$  (the closure  of $C$ in $Z$), so $\bar{C}\setminus C={\rm cl}_Z(C)\setminus C\neq \emptyset $  and the result follows by the assumption on $X.$

Let $U$ be an open  definable subset of $X$ and let $C\subseteq U$ be a definable subset which is not closed in $U$.  Note that  $\bar{C}\cap U={\rm cl}_U(C)$  (the closure of $C$ in $U$). Let $B=C\cup (\bar{C}\setminus U)=C\cup ((\bar{C}\setminus C)\setminus U).$ Then we have $\emptyset \neq {\rm cl}_U(C)\setminus C=\bar{C}\cap U\setminus C =(\bar{C}\setminus C)\cap U= \bar{B}\setminus B$  and the result follows applying  the assumption on $X$ to $B.$ Note that  any definable curve in $B$ with limit a point in $\bar{B}\setminus B\subseteq U$ must enter $U$ and so gives a definable curve in $C=B\cap U.$

Let $Z\cap U$ be a general  locally closed definable subset  of $X$, where $Z$ is a closed definable subset and $U$ is an open definable subset. Let $C\subseteq Z\cap U$ be a definable subset which is not closed in $Z\cap U.$ Then ${\rm cl}_{Z\cap U}(C)=\bar{C}\cap U={\rm cl}_U(C)$ and  ${\rm cl}_U(C)\setminus C={\rm cl}_{Z\cap U}(C)\setminus C\neq \emptyset $ and therefore,  the result follows from the previous case.
\qed \\

\begin{lem}\label{lem aecs union opens}
Suppose that $X$ is a  definable space and $V$ and $W$ are open definable subsets such that $V\cup W=X$ and almost everywhere curve selection holds for $V$ and $W$. Then almost everywhere curve selection holds for $X$.
\end{lem}

\pf
Let $C\subseteq X=V\cup W$ be a definable subset  which is not closed. Let $C_V=C\cap V\subseteq V$ and let $C_W=C\cap W\subseteq W$. Then we have $C=C_V\cup C_W$, $\bar{C}=\bar{C_V}\cup \bar{C_W}$ and  ${\rm cl}_V(C_V)=\bar{C_V}\cap V=\bar{C}\cap V$  and similarly  ${\rm cl}_W(C_W)=\bar{C_W}\cap W=\bar{C}\cap W.$ 
So  $\bar{C}=(\bar{C}\cap V)\cup (\bar{C}\cap W)={\rm cl}_V(C_V) \cup {\rm cl}_W(C_W).$ Therefore, $\bar{C}\setminus C=({\rm cl}_V(C_V)\setminus C) \cup ({\rm cl}_W(C_W)\setminus C)=({\rm cl}_V(C_V)\setminus C_V) \cup ({\rm cl}_W(C_W)\setminus C_W).$ 

If $C_V$ is not closed in $V$, by the hypothesis, there is a definable set $F_V\subseteq {\rm cl}_V(C_V)\setminus C_V$ such that  $\dim F_V<\dim({\rm cl}_V(C_V)\setminus C_V)$ and  for every $x\in  {\rm cl}_V(C_V)\setminus (C_V\cup F_V)$ there is a definable curve in $C_V$ which has $x$ as a limit point. Similarly, if $C_W$ is not closed in $W$, there is a definable set $F_W\subseteq {\rm cl}_W(C_W)\setminus C_W$ such that  $\dim F_W<\dim({\rm cl}_W(C_W)\setminus C_W)$ and  for every $x\in  {\rm cl}_W(C_W)\setminus (C_W\cup F_W)$ there is a definable curve in $C_W$ which has $x$ as a limit point. Let $E_V$ be $F_V$ if it exists and let it be $\emptyset $ otherwise. Similarly, let $E_W$ be $F_W$ if it exists and let it be $\emptyset $ otherwise. Let $E=E_V\cup E_W$. Since $\bar{C}\setminus C=({\rm cl}_V(C_V)\setminus C_V) \cup ({\rm cl}_W(C_W)\setminus C_W)$ we have  $E\subseteq \bar{C}\setminus C.$ Since $C=C_V\cup C_W$ we also have that for every $x\in  \bar{C}\setminus (C\cup E)$ there is a definable curve in $C$ which has $x$ as a limit point. Since $\dim E=\max \{\dim E_V, \dim E_W\}$ and $\dim \bar{C}\setminus C=\max \{\dim({\rm cl}_V(C_V)\setminus C_V), \dim({\rm cl}_W(C)\setminus C_W)\}$ we also have $\dim E<\dim(\bar{C}\setminus C)$ as required. \qed \\

By Fact \ref{fact arcs ps}, Lemma \ref{lem aecs loc closed} and an induction argument using  Lemma \ref{lem aecs union opens} we see that:

\begin{cor}\label{cor aecs man}
Almost everywhere curve selection holds for locally closed definable subsets of definable manifolds. \\
\end{cor}

Let $X$ be a definable space. We say that:
\begin{itemize}
\item
 $X$ is {\it locally definably compact} if every $x\in X$ has  a definably compact neighborhood.\\
\end{itemize}


We now have the following  extension of the almost everywhere curve selection to the non-affine case which will be useful later:

\begin{thm}[Local almost everywhere curve selection]\label{thm local aecs}
Suppose that ${\mathbb M}$ has definable Skolem functions. Let $X$ be a  Hausdorff, locally definably compact definable space. If $C\subseteq X$ is a definable subset which is not closed, then for every $z\in \bar{C}\setminus C$ there is a definable open neighborhood $V$ of $z$ in $X$ such that $\bar{V}$ is definably compact and there is a definable set $E\subseteq (\bar{C}\setminus C)\cap V$ such that  $\dim E<\dim((\bar{C}\cap V)\setminus (C\cap V))$ and  for every $x\in  (\bar{C}\cap V)\setminus ((C\cap V)\cup E)$ there is a definable curve in $C\cap V$ which has $x$ as a limit point. 
\end{thm}

\pf
By the assumption on $X$ we get $V$ such that $\bar{V}$ is definably compact and so definably normal (Theorem \ref{thm def comp hausd normal}). The result  then follows at once  after we show that almost everywhere curve selection holds for definably normal, definably compact definable spaces $Y.$ 

So let $Y$ be such a definable space. Consider the definable charts $(U_i,\phi _i)_{i=1}^l$ of $Y$. Since $Y$ is definably normal, by the shrinking lemma, there are open definable subsets $V_i$ ($1\leq i \leq l$) and closed definable subsets $C_i$ ($1\leq i \leq l$) such that $V_i\subseteq C_i\subseteq U_i$ and $Y=\cup \{V_i:i=1, \ldots , l\}.$ Since each $C_i$ is definably compact and each $\phi _i$ is a definable homeomorphism,  we have that each $\phi _i(C_i)$ is a closed  (and bounded) definable subset of $M^{n_i}$ and so  by Fact \ref{fact arcs ps}  and Lemma \ref{lem aecs loc closed}, each $\phi _i(C_i)$ and so each $C_i$ has almost everywhere curve selection. So by Lemma \ref{lem aecs loc closed}, each $V_i$ has almost everywhere curve selection. Now we conclude by induction of $l,$ using Lemma \ref{lem aecs union opens}, that $Y$ has almost everywhere curve selection. 
\qed \\

The second part of the proof of Theorem \ref{thm local aecs} shows that: 

\begin{cor}\label{cor aecs comp normal}
Almost everywhere curve selection holds for  definably normal,  definably compact definable spaces - even without assuming that ${\mathbb M}$ has definable Skolem functions.\\
\end{cor}

\end{subsection}

\end{section}

\begin{section}{Proper morphisms in $\Df$}\label{section  def proper maps}

\begin{subsection}{Preliminaries}\label{section prelim}
Here we recall some preliminary notions for the category $\Df$ whose objects are definable spaces and whose morphism are continuous definable maps between definable spaces.\\

Let $f:X\to Y$ be a morphism in $\Df.$   We say that:
\begin{itemize}
\item
$f:X\to Y$ is  {\it closed in $\Df$ (i.e., definably closed)} if for every object $A$ of $\Df$ such that $A$ is a closed  subset  of $X$, its image $f(A)$ is a closed (definable) subset of $Y.$ 
\item
$f:X\to Y$ is  a {\it closed (resp. open)  immersion} if $f:X\to f(X)$ is a  homeomorphism and $f(X)$ is  a closed (resp. open) subset of $Y$.\\
\end{itemize}

\begin{prop}
In the category $\Df$  the {\it cartesian square} of any two morphisms $f:X\to Z$ and $g:Y\to Z$ in $\Df$ exists and is given by a commutative diagram
$$
\xymatrix{X\times _ZY \ar[r]^{p_Y} \ar[d]^{p_X} & Y \ar[d]^g \\
X \ar[r]^f &Z}
$$
where the morphisms $p_X$ and $p_Y$ are known as projections. The cartesian square satisfies the following universal property:   for any other object $Q$ of $\Df$ and morphisms $q_X:Q\to X$ and $q_Y:Q\to Y$ of $\Df$ for which the following diagram commutes,
$$
\xymatrix{
Q \ar@/_/[ddr]_{q_X} \ar@/^/[drr]^{q_Y}
\ar@{-->}[dr]^-{u} \\
& X\times _ZY \ar[d]^{p_X} \ar[r]^{p_Y}
& Y \ar[d]^g \\
& X \ar[r]^f & Z }
$$
there  exist a unique natural morphism $u:Q\to X\times _ZY$ (called mediating morphism) making the whole diagram commute. As with all universal constructions, the cartesian square  is unique up to a definable homeomorphism.\\
\end{prop}

\pf
The usual fiber product $X\times _ZY=\{\langle x,y\rangle \in X\times Y: f(x)=g(y)\}$ (a closed  definable subspace of the definable space $X\times Y$) together with the restrictions  $p_X:X\times _ZY\to X$ and $p_Y:X\times _ZY\to Y$ of the usual projections determine a cartesian square in the category $\Df$. 
\qed \\



Given a morphism $f:X\to Y$ in $\Df$, the corresponding diagonal morphism is the unique morphism $\Delta :X\to X\times _YX$ in $\Df$ given by the universal property of cartesian squares:
$$
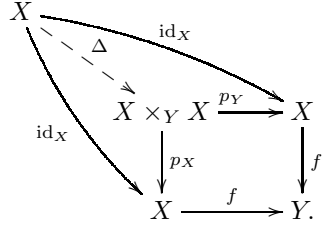
\xymatrix{
X \ar@/_/[ddr]_{{\rm id}_X} \ar@/^/[drr]^{{\rm id}_X}
\ar@{-->}[dr]^-{\Delta } \\
& X\times _YX \ar[d]^{p_X} \ar[r]^{p_Y}
& X \ar[d]^f \\
& X \ar[r]^f & Y. }
$$
We say that:
\begin{itemize}
\item
 $f:X\to Y$ is {\it separated in $\Df$} if the corresponding diagonal morphism  $\Delta :X\to X\times _YX$ is a closed  immersion. \\
 \end{itemize}
 
 \noindent
 We say that an object $Z$ in $\Df$ is {\it separated in $\Df$} if the morphism $Z\to \{{\rm pt}\}$ to a point is separated.

\begin{nrmk}\label{nrmk sep df}
{\em
Since in the above diagram we have $p_X\circ \Delta =p_Y\circ \Delta ={\rm id}_X$, it is clear that the following are equivalent:
\begin{enumerate}

\item
$f:X\to Y$ (resp. $Z$)  is separated in $\Df$.

\item
The image of the corresponding diagonal morphism  $\Delta :X\to X\times _YX$ is a  closed  (definable) subset of $X\times _YX$ (resp. the diagonal $\Delta _Z$ of $Z$ is a  closed  (definable) subset of $Z\times Z$).\\
\end{enumerate}
}
\end{nrmk}

Let  $Z$ be an object of $\Df$ and $s:Z'\to Z$  a morphism in $\Df.$
\begin{itemize}
\item
 By   a {\it morphism over $Z$ in $\Df$} we mean  a commutative diagram
$$
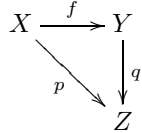
\xymatrix{ X \ar [rd]_p \ar[r]^f & Y \ar[d]^q \\
& Z}
$$
of morphisms in $\Df.$
\item
We call $s:Z'\to Z$   a {\it base extension in $\Df$} and the induced commutative diagram
$$
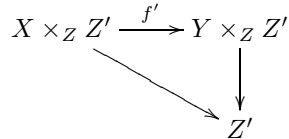
\xymatrix{ X\times _ZZ' \ar [rd] \ar[r]^{f'} & Y\times _ZZ' \ar[d] \\
& Z'}
$$
where $f'=f\times {\rm id}_{Z'} $ and the downarrows are the natural projections,  is called  the {\it base extension in $\Df$ of} 
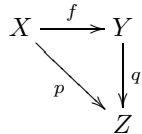
$$
\xymatrix{ X \ar [rd]_p \ar[r]^f & Y \ar[d]^q \\
& Z}
$$
Note that since the $f'$ above is completely determined by the corresponding morphism over $Z$ we will often just say that $f':X\times _ZZ'\to Y\times _ZZ'$ is the corresponding base extension  morphism.\\
\end{itemize}

Let $f:X\to Y$ be a morphism in $\Df.$   We say that:
\begin{itemize}
\item 
$f:X\to Y$ is {\it universally  closed in $\Df$} if for any morphism $g:Y'\to Y$  in $\Df$ the morphism $f':X'\to Y'$ in $\Df$ obtained from the cartesian square 
$$
\xymatrix{X' \ar[r]^{f'} \ar[d]^{g'} & Y' \ar[d]^g \\
X \ar[r]^f &Y}
$$
in $\Df$ is   closed in $\Df$.  \\
\end{itemize}

\begin{defn}\label{defn proper in df}
{\em
We say that a morphism $f:X\to Y$ in $\Df$ is {\it proper in $\Df$} if $f:X\to Y$ is separated and universally closed  in $\Df$.\\
}
\end{defn}

\begin{defn}\label{defn complete in df}
{\em
We say that an object $Z$ of $\Df$ is {\it complete in $\Df$} if the morphism $Z\to {\rm pt}$  is proper in $\Df.$ \\
}
\end{defn}

Below we will relate the notion of proper in $\Df$ and complete in $\Df$ with the usual notions of definably proper and definably compact.

\end{subsection}

\begin{subsection}{Separated and proper in $\Df$}\label{subsection sep and proper in df} Here we list the main properties of morphisms  separated or proper in $\Df$.\\

From Remark \ref{nrmk sep df} and the way cartesian squares are defined in $\Df$ we easily obtain the following:\\

\begin{nrmk}\label{nrmk sep df2}
{\em
Let $f:X\to Y$ be a morphism in $\Df$. 
Then the  following are equivalent:
\begin{enumerate}

\item
$f:X\to Y$  is separated in $\Df$.

\item
The fibers $f^{-1}(y)$ of $f$ are Hausdorff (with the induced topology). \\
\end{enumerate}
}
\end{nrmk}

 Directly from the definitions (as in \cite[Chapter I, Propositions 5.5.1 and 5.5.5]{ega1}) or more easily from Remark \ref{nrmk sep df2}   the following  is immediate:\\

\begin{prop}\label{prop def separated4} 
In the category $\Df$ the following hold:
\begin{enumerate}
\item
Injective continuous definable maps are separated in $\Df$.
\item
A composition of two morphisms separated in $\Df$   is separated in $\Df$.
\item
Let $
\xymatrix{ X \ar [rd]_p \ar[r]^f & Y \ar[d]^q \\
& Z}
$
be a morphism over $Z$ in $\Df$ and $Z'\to Z$ a base extension in $\Df.$ 
If $f:X\to Y$ is  separated   in $\Df,$  then the corresponding base extension morphism $f':X\times _ZZ'\to Y\times _ZZ'$ is separated in $\Df$.
\item
Let $
\xymatrix{ X \ar [rd]_p \ar[r]^f & Y \ar[d]^q \\
& Z}
$
and 
$
\xymatrix{ X' \ar [rd]_{p'} \ar[r]^{f'} & Y' \ar[d]^{q'} \\
& Z}
$
be morphisms over $Z$ in $\Df.$ 
If $f:X\to Y$ and $f':X'\to Y'$ are  separated in $\Df$, then the corresponding product morphism $f\times f':X\times _ZX' \to Y\times _ZY'$   is separated in $\Df$.
\item
If $f:X\to Y$ and $g:Y\to Z$ are morphisms such that $g\circ f$ is separated in $\Df$, then $f$ is separated in $\Df$.
\item
A morphism $f:X\to Y$  is separated in $\Df$ if and only if $Y$ can be covered by finitely many open definable subsets $V_i$ such that $f_|:f^{-1}(V_i)\to V_i$ is separated in $\Df$.\\
\end{enumerate}
\end{prop}

Directly from the definitions (as in \cite[Chapter II, Proposition 5.4.2 and Corollary 5.4.3]{ega2}, see also \cite[Section 9]{dk1}) one has   the following. For the readers convenience we include some details:\\

\begin{prop}\label{prop proper in def} 
In the category $\Df$ the following hold:
\begin{enumerate}
\item
Closed  immersions are proper in $\Df$.
\item
A composition of two morphisms  proper in $\Df$   is proper in $\Df$.
\item
Let $
\xymatrix{ X \ar [rd]_p \ar[r]^f & Y \ar[d]^q \\
& Z}
$
be a morphism over $Z$ in $\Df$ and $Z'\to Z$ a base extension in $\Df.$ 
If $f:X\to Y$ is  proper  in $\Df,$  then the corresponding base extension morphism $f':X\times _ZZ'\to Y\times _ZZ'$ is proper in $\Df$.
\item
Let $
\xymatrix{ X \ar [rd]_p \ar[r]^f & Y \ar[d]^q \\
& Z}
$
and 
$
\xymatrix{ X' \ar [rd]_{p'} \ar[r]^{f'} & Y' \ar[d]^{q'} \\
& Z}
$
be morphisms over $Z$ in $\Df.$ 
If $f:X\to Y$ and $f':X'\to Y'$ are proper in $\Df$, then the corresponding product morphism $f\times f':X\times _ZX' \to Y\times _ZY'$   is proper in $\Df$.
\item
If $f:X\to Y$ and $g:Y\to Z$ are morphisms such that $g\circ f$ is proper in $\Df$, then:
\begin{itemize}
\item[(i)]
$f$ is proper in $\Df;$
 \item[(ii)]
 if $g$ is separated in $\Df$ and $f$ is surjective, then $g$ proper in $\Df.$
 \end{itemize}
\item
A morphism $f:X\to Y$  is proper in $\Df$ if and only if $Y$ can be covered by finitely many open definable subsets $V_i$ such that $f_|:f^{-1}(V_i)\to V_i$ is proper in $\Df$.\\
\end{enumerate}
\end{prop}

\pf
(1) Let $X\to Y$ be a closed immersion and $Y'\to Y$   a  morphism in $\Df.$ Since $X\times _Y Y'\to Y\times _Y Y'=Y'$ is also a closed immersion, it is closed in $\Df.$ So $X\to Y$ is universally closed and is separated by Proposition \ref{prop def separated4} (1).

(2) Let $X\to Y$ and $Y\to Z$ be morphisms proper in $\Df$ and let $Z'\to Z$ be a morphism in $\Df.$ Since $X\times _{Z}Z'=X\times _Y (Y\times _ZZ')$ and $X\times _ZZ'\to Z'$ is $X\times _Y(Y\times _Z Z')\to Y\times _Z Z' \to Z'$, 
the result follows from the fact that the composition of morphisms closed in $\Df$ is closed in $\Df$ and Proposition \ref{prop def separated4} (2).

(3) Let $
\xymatrix{ X \ar [rd] \ar[r] & Y \ar[d] \\
& Z}
$
be a morphism over $Z$ in $\Df,$  $Z'\to Z$ a base extension in $\Df$ and suppose that    $X\to Y$ is a morphism proper in $\Df.$  Since $X\times _ZZ'= X\times _Y(Y\times _ZZ'),$ for every morphism $W\to Y\times _ZZ'$ we have 
$$(X\times _ZZ')\times _{Y\times _ZZ'}Z = (X\times _Y(Y\times _ZZ'))\times _{Y\times _ZZ'}W=X\times _YW.$$
Hence, since $X\times _YW\to W$ is closed in $\Df$ by hypothesis, the result follows using also Proposition \ref{prop def separated4} (3).

(4) Let $
\xymatrix{ X \ar [rd] \ar[r] & Y \ar[d] \\
& Z}
$
and 
$
\xymatrix{ X' \ar [rd]_{} \ar[r]^{} & Y' \ar[d]^{} \\
& Z}
$
be morphisms over $Z$ in $\Df$ with $X\to Y$ and $X'\to Y'$  proper in $\Df$. Then the product morphism $X\times _ZX'\to Y\times _ZY'$ is the composition of the base extension $X\times _ZX'\to Y\times _ZX',$ the identification $Y\times _ZX'=X'\times _ZY$ and the base extension $X'\times _ZY\to Y'\times _ZY.$ So the result follows from (1) and (3).

(5) Let $X\to Y$ and $Y\to Z$ be morphisms in $\Df$ such that the composition $X\to Y\to Z$ is proper in $\Df$. 

(i) Let $Y'\to Y$ be a morphism in $\Df$. Then $X\times _ZY' \to Y'$ obtained with the composition $Y'\to Y\to Z$ is the same as $X\times _YY'\to Y'$. So since $X\times _ZY'\to Y'$ is closed in $\Df$, so is $X\times _YY'\to Y'$ and the result follows using also Proposition \ref{prop def separated4} (5). 

(ii) Let $Z'\to Z$ be a morphism in $\Df$. Then 
$$
\xymatrix{X\times _ZZ' \ar[r]^{f\times {\rm id}_{Z'}} \ar[dr]^{p} & Y\times _ZZ' \ar[d]^{p'} \\
  & Z'}
$$
is a commutative diagram, with $f\times {\rm id}_{Z'}$ surjective and $p$ closed in $\Df$ by hypothesis. It follows that $p'$ is closed in $\Df$ as required.

(6) Suppose that $f:X\to Y$ is a morphism in $\Df$ and let $\{V_i\}_{i\leq k}$ be a finite cover of $Y$ by open definable subsets. If $g:Y'\to Y$ is a morphism in $\Df$, then $\{f^{-1}(V_i)\}_{i\leq k}$ (resp. $\{g^{-1}(V_i)\}_{i\leq k}$) is a finite cover of $X$ (resp. $Y'$) by open definable subsets and  $\{f^{-1}(V_i)\times _{Y}g^{-1}(V_i)\}_{i\leq k}$ is a finite cover of $X\times _YY'$ by open definable subsets. One the other hand, $f^{-1}(V_i)\times _{Y}g^{-1}(V_i)= f^{-1}(V_i)\times _{V_i} g^{-1}(V_i)$ and 
$$
\xymatrix{f^{-1}(V_i)\times _{V_i}g^{-1}(V_i) \ar[r]^{i} \ar[d]^{p'_i} & X\times _YY' \ar[d]^{p'} \\
g^{-1}(V_i)  \ar[r]^{j}& Y'}
$$
is a commutative diagram with  $i$ and $j$ the inclusions, $p'$ the projection and $p'_i$ the restriction of $p'$.  Since $p'$ is closed in $\Df$ if and only if each $p'_i$ is closed in $\Df$ the result follows using also Proposition \ref{prop def separated4} (6).
\qed \\

\begin{cor}\label{cor basic compl}
Let $f:X\to Y$ be a morphism in $\Df$ and $Z\subseteq X$ an object in $\Df$ which is  complete in $\Df$. Then the following hold:
\begin{enumerate}
\item 
$Z$ is a closed (definable) subset of $X.$
\item
$f_{|Z}:Z\to Y$ is proper in $\Df.$
\item
$f(Z)\subseteq Y$ is  (definable) complete in $\Df.$ 
\item
If $f:X\to Y$ is proper in $\Df$ and $C\subseteq Y$ is an object in $\Df$ which is  complete in $\Df,$ then $f^{-1}(C)\subseteq X$  is  (definable) complete in $\Df.$ \\
\end{enumerate}
\end{cor}

From Proposition \ref{prop proper in def} we also obtain in a standard way the following:\\

\begin{cor}\label{cor completion in B}
Let $\bB$ be a full a subcategory  of the category of definable spaces $\Df$ whose set of objects is:
\begin{itemize}
\item
closed under taking locally closed definable subspaces of objects of $\bB$,
\item
closed under taking cartesian products of objects of $\bB.$
\end{itemize}
Then the following are equivalent:
\begin{enumerate}
\item
Every object $X$ of $\bB$ is {\it completable in $\bB$}  i.e., there exists an object $X'$ of $\bB$ which is complete in $\Df$ together with an open immersion $i: X\hookrightarrow X'$ in $\bB$ with  $i(X)$ dense in $X'$. Such $i: X\hookrightarrow X'$ is called a {\it completion} of $X$ in $\bB$.

\item
Every  morphism $f:X\to Y$ in $\bB$ is {\it  completable in $\bB$} i..e,  there exists a commutative diagram
$$
\xymatrix{ X \ar [d]_f \ar[r]^-{i} & X' \ar[d]^{f'} \\
Y \ar[r]^-{j} & Y'}
$$
of  morphisms in $\bB$ such that:  (i) $i:X\to X'$ is a completion of $X$  in $\bB$; (ii) $j$ is a  completion of $Y$ in $\bB$.
\item
Every morphism $f:X\to Y$ in $\bB$ has {\it a proper extension in $\bB$} i.e.,  there exists a commutative diagram
$$
\xymatrix{ X \ar [rd]_f \ar[r]^\iota & P \ar[d]^{\overline{f}} \\
& Y}
$$
of morphisms in $\bB$ such that $\iota$ is a  open immersion with $\iota (X)$ dense in $P$ and $\overline{f}$ a  proper  in $\Df$.\\
\end{enumerate}
\end{cor}

\pf
Assume (1). Let $h:X\to Y$ be a morphism in ${\mathbf B}$. Let $j:Y\to Y'$ be a  completion of $Y$ in ${\mathbf B}$. Choose also a completion $g:X\to X''$ of $X$ in ${\mathbf B}$ and note that $g\times j: X\times Y\to X''\times Y'$ is a  completion of $X\times Y$ in ${\mathbf B}$ (since $X''\times Y'$ is complete in $\Df$ by Proposition \ref{prop proper in def}  (4)). 
Let $X'$ be the closure of $(g\times j)(\Gamma (h))$ in $X''\times Y'.$ Then $i:X\to X'$ given by $i=(g\times j)\circ ({\rm id}_X\times h)$ is a  completion of $X$ in ${\mathbf B}$ (by Proposition \ref{prop proper in def}  (1) and (5)), and the restriction of the projection $X''\times Y'\to Y'$ to $X'$ is a morphism $h':X' \to Y'$ making a commutative diagram
$$
\xymatrix{ X \ar [d]_h \ar[r]^-{i} & X' \ar[d]^{h'} \\
Y \ar[r]^-{j} & Y'}
$$
of  morphisms in ${\mathbf B}$ as required in (2). 

Assume (2). Let $h:X\to Y$ be a morphism in ${\mathbf B}$.  Then there exists a commutative diagram
$$
\xymatrix{ X \ar [d]_h \ar[r]^-{i} & X' \ar[d]^{h'} \\
Y \ar[r]^-{j} & Y'}
$$
of   morphisms in ${\mathbf B}$ such that: (i) $i:X\to X'$ is a  completion of $X$  in ${\mathbf B}$; (ii) $j:Y\to Y'$ is a   completion of $Y$ in ${\mathbf B}$.  Let $P=h'^{-1}(j(Y))$ (an open definable subspace of $X'$) and $\bar{h}=j^{-1}\circ h'_{|P}:P\to Y$ where $j^{-1}:j(Y)\to Y$ is the inverse of $j:Y\to j(Y)$ which is a definable homeomorphism. Then we have a commutative diagram
$$
\xymatrix{ X \ar [rd]_h \ar[r]^\iota & P \ar[d]^{\overline{h}} \\
& Y}
$$
of  morphisms in ${\mathbf B}$ such that $\iota =i:X\to P$ is a definable  open immersion with $\iota (X)$ dense in $P$ and $\overline{h}$ is  proper in $\Df$ (since $h':X'\to Y'$ is proper in $\Df$ by Corollary \ref{cor basic compl} (2)) as required in (3).

Assume (3). Let $X$ an object of ${\mathbf B}$. Take $h:X\to \{{\rm pt}\}$ to be the morphism in ${\mathbf B}$  to a point. Applying (3) to this morphism we obtain (1).
\qed \\

\end{subsection}

\begin{subsection}{Definably proper maps}\label{subsection def proper}
Here we recall the definition of definably proper map between definable spaces and prove its main properties. A special case of this theory appears in \cite[Chapter 6, Section 4]{vdd} in the context of affine definable spaces in o-minimal expansions of ordered groups. \\

\begin{defn}\label{defn def proper}
{\em
A continuos definable map $f:X\to Y$ between definable spaces $X$ and $Y$ is called {\it definably proper} if for every definably compact definable subset $K$ of $Y$ its inverse image $f^{-1}(K)$ is a definably compact definable subset of $X$.
}
\end{defn}

From the definitions we see that:

\begin{nrmk}\label{nrmk proper compact}
A definable space $X$ is definably compact if and only if the map $X\to \{\pt\}$ to a point is definably proper.
\end{nrmk}

Typical examples of  definably proper continuous definable maps are: (i) $f:X\to Y$ where $X$ is a definably compact definable space and $Y$ is any definable space; (ii)  the projection $X\times Y\to Y$ where $X$ is a definably compact definable space and $Y$ is any definable space; (iii) closed definable immersions.\\

The following is proved just like in the affine case in o-minimal expansions of ordered groups treated in \cite[Chapter 6, Lemma (4.5)]{vdd}:\\

\begin{thm}\label{thm def proper and curves}
Let $f:X\to Y$ be a continuous definable map. Suppose that every definably compact subset of $Y$ is a closed subset (e.g. ${\mathbb M}$ has definable Skolem functions and $Y$ is Hausdorff). Then the following are equivalent:
\begin{enumerate}
\item
 $f$ is definably proper.

 \item
For every definable curve $\alpha :(a,b)\to  X$ and every continuous definable map $[a,b]\to Y$ which makes a commutative diagram
$$
\xymatrix{
(a,b)  \ar[r]^-{\alpha } \ar@{^{(}->}[d] &  X \ar[d]^-f\\
[a,b] \ar[r] \ar@{-->}[ur]_-{\overline{\alpha }} & Y }
$$
there is at least one continuous definable map $[a,b]\to X$ making the whole diagram commutative.
\end{enumerate}
\end{thm}


\pf
Assume (1). Let $\alpha : (a,b)\to X$ be a definable curve in $X$ such that $f\circ \alpha  $ is completable in $Y$, say $\lim _{t\to b^-}f\circ \alpha  (t)=y\in Y$.  Take $c\in (a,b)$ and set $K=\{f(\alpha  (t)):t\in [c,b)\}\cup \{y\}\subseteq Y.$ Then $K$ is a definably compact definable subset of $Y$ and so, $f^{-1}(K)$ is a definably compact definable subset of $X$ containing $\alpha((c,b))$. Thus $\alpha  $ must be completable in $f^{-1}(K)$, hence in $X$.

Assume (2). Suppose that $f$ is not definably proper. Then there is a definably compact definable subset $K$ of $Y$ such that $f^{-1}(K)$ is not a definably compact definable subset of $X$. Thus there is a definable curve $\alpha  :(a,b)\to f^{-1}(K)\subseteq X$ in $f^{-1}(K)$ which is not completable in $f^{-1}(K)$. Since $f^{-1}(K)$ is closed (by assumption on $Y$, $K$ is closed), $\alpha  $ is not completable. But $f\circ \alpha :(a,b)\to K\subseteq Y$ is completable which contradicts (2).  \qed \\

By Theorem \ref{thm def proper and curves} we  have the following which  summarizes the main properties of definably proper maps.

\begin{cor}\label{cor def proper4} 
Let ${\mathbf A}$ be a  full subcategory of $\Df$ such that every definably compact subset of an object of ${\mathbf A}$ is a closed subset. Suppose that the set of objects of ${\mathbf A}$ is:
 \begin{itemize}
 \item
 closed under taking  locally closed definable subsets of objects of ${\mathbf A}$;
 \item
 is closed under taking cartesian products of objects of ${\mathbf A}$.
 \end{itemize}
 In the category ${\mathbf A}$ the following hold:
\begin{enumerate}
\item
Closed  immersions are definably proper.
\item
A composition of two definably proper morphisms  is definably proper.
\item
Let $
\xymatrix{ X \ar [rd]_p \ar[r]^f & Y \ar[d]^q \\
& Z}
$
be a morphism over $Z$ in ${\mathbf A}$ and $Z'\to Z$ a base extension in ${\mathbf A}.$ 
If $f:X\to Y$ is   definably proper, then the corresponding base extension morphism $f':X\times _ZZ'\to Y\times _ZZ'$ is definably proper.
\item
Let $
\xymatrix{ X \ar [rd]_p \ar[r]^f & Y \ar[d]^q \\
& Z}
$
and 
$
\xymatrix{ X' \ar [rd]_{p'} \ar[r]^{f'} & Y' \ar[d]^{q'} \\
& Z}
$
be morphisms over $Z$ in ${\mathbf A}.$ 
If $f:X\to Y$ and $f':X'\to Y'$ are definably proper, then the corresponding product morphism $f\times f':X\times _ZX' \to Y\times _ZY'$   is definably proper.
\item
If $f:X\to Y$ and $g:Y\to Z$ are morphisms such that $g\circ f$ is definably proper, then:
\begin{itemize} 
\item[(i)]
$f$ is definably proper;
\item[(ii)]
if ${\mathbb M}$ has definable Skolem functions, then $g_{|f(X)}:f(X)\to Z$ is definably proper.
\end{itemize}
\item
A morphism $f:X\to Y$  is definably proper if and only if $Y$ can be covered by finitely many open definable subsets $V_i$ such that $f_|:f^{-1}(V_i)\to V_i$ is definably proper.
\end{enumerate}
\end{cor}

\pf
(1) Consider the commutative diagram:
$$
\xymatrix{
(a,b) \ar[dd]^-{} \ar[r]^{\alpha } & X \ar[d]^-{f}\\
& f(X)\ar[d]^-{} \\
[a,b] \ar@{-->}[uur]^-{\gamma '} \ar@{-->}[ur]^{\gamma }\ar[r]^{\bar{\alpha }} & Y 
}
$$
where $f:X\to Y$ is a definable closed immersion and we assume we have $\alpha $ such that $\bar{\alpha}$ exists. We must show that $\gamma '$ exists. Since  the inclusion $f(X)\subseteq Y$  is closed  and we have $f\circ \alpha $ such that $\bar{\alpha }$ exists,   $\gamma $ exists. So, since $f:X\to f(X)$ is a definable homeomorphism, we let $\gamma '=f^{-1}\circ \gamma $.

(2) Consider the commutative diagram:
$$
\xymatrix{
(a,b) \ar[dd]^-{} \ar[r]^{\alpha } & X \ar[d]^-{f}\\
& Y\ar[d]^-{g} \\
[a,b] \ar@{-->}[uur]^-{\gamma '} \ar@{-->}[ur]^{\gamma }\ar[r]^{\bar{\alpha }} & Z
}
$$
where we assume we have $\alpha $ such that $\bar{\alpha}$ exists. We must show that $\gamma '$ exists. Since  $g:Y\to Z$ is definably proper   and we have $f\circ \alpha $ such that $\bar{\alpha }$ exists,  by Theorem \ref{thm def proper and curves} $\gamma $ exists. Since  $f:X\to Y$ is definably proper   and we have $\alpha $ such that $\gamma $ exists,  by Theorem \ref{thm def proper and curves} $\gamma '$ exists.

(3) Since the base extension morphism is a special case of the product morphism, the result follows from (4) below.

(4) Consider the commutative diagram:
$$
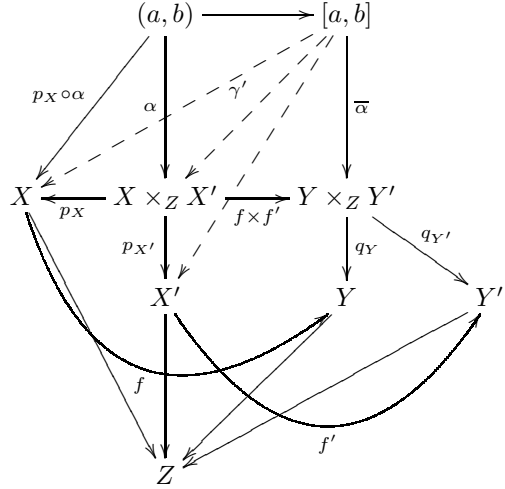
\xymatrix{
& \ar[ddl]_-{p_X\circ \alpha }(a,b) \ar[dd]_-{\alpha } \ar[r]^{} &  \ar@{-->}[ddll]_-{}  \ar@{-->}[ddl]_-{\gamma '} \ar@{-->}[dddl]_-{}[a,b] \ar[dd]^-{\bar{\alpha }} & \\
& & & \\
X  \ar[dddr]^-{} \ar@/_4.0pc/[drr]_-{f} & \ar[l]^-{p_X}X\times _ZX' \ar[d]_-{p_{X'}}\ar[r]_-{f\times f'}& Y\times _ZY'\ar[d]^-{q_Y}\ar[dr]^-{q_{Y'}} & \\
 & X' \ar[dd]^-{}\ar@/_4.0pc/[rr]_-{f'}& \ar[ddl]^-{}Y & \ar[ddll]^-{}Y' \\
& & & \\
& Z & & 
}
$$
where we assume we have $\alpha $ such that $\bar{\alpha}$ exists. We must show that $\gamma ': [a, b]\to X\times _ZX'$ exists. Since  $f:X\to Y$ is definably proper   and we have $p_X\circ \alpha $ such that $q_Y\circ \bar{\alpha }$ exists,  by Theorem \ref{thm def proper and curves}, $[a,b] \to X$ exists. Since  $f':X'\to Y'$ is definably proper   and we have $p_{X'}\circ \alpha $ such that $q_{Y'}\circ \bar{\alpha} $ exists,  by Theorem \ref{thm def proper and curves} $[a,b]\to X'$ exists. So we let $\gamma '$ be the morphism given by the universal property of Cartesian squares.

(5) 

(i) Consider the commutative diagram:
$$
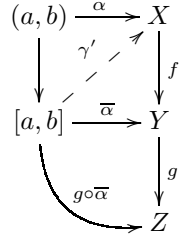
\xymatrix{
(a,b) \ar[d]^-{} \ar[r]^{\alpha } & X \ar[d]^f \\
[a,b] \ar@/_1.5pc/[dr]^-{g\circ \bar{\alpha}}\ar@{-->}[ur]^-{\gamma '} \ar[r]^-{\bar{\alpha} } & Y \ar[d]^g  \\
& Z
}
$$
where we assume we have $\alpha $ such that $\bar{\alpha}$ exists. We must show that $\gamma '$ exists. Since  $g\circ f:X\to Y$  is definably proper and we have $\alpha $ such that $g\circ \bar{\alpha }$ exists, by Theorem \ref{thm def proper and curves},  $\gamma '$ exists. 

(ii) Consider the commutative diagram:
$$
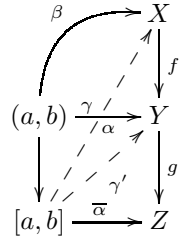
\xymatrix{
& X \ar[d]^f \\
(a,b) \ar[d] \ar@/^1.5pc/[ur]^-{\beta  } \ar[r]_-{\alpha } & Y \ar[d]^g  \\
[a,b] \ar@{-->}[uur]^-{\gamma }\ar@{-->}[ur]_-{\gamma '} \ar[r]^{\bar{\alpha }}& Z
}
$$
where we assume we have $\alpha $ such that $\bar{\alpha}$ exists. We must show that $\gamma '$ exists. Since $f$ is surjective, by definable Skolem function let $\beta $ be such that $\alpha =f\circ \beta .$ Since  $g\circ f:X\to Y$  is definably proper and we have $\beta  $ such that $ \bar{\alpha }$ exists, by Theorem \ref{thm def proper and curves},  $\gamma $ exists. Now take $\gamma '=f\circ \gamma .$

(6) One implication is clear. Suppose that there are open definable subsets $V_1, \ldots , V_l$ of $Y$ such that each restriction $f_{|}:f^{-1}(V_i)\to V_i$ is definably proper. Let $\alpha : (a,b)\to X$ be a definable curve such that $f\circ \alpha :(a,b)\to Y$ is completable. Without loss of generality it is enough to  show that $\lim_{t\to b^-}\alpha (t)$ exists in $X$. Let $z=\lim _{t\to b^-}f\circ \alpha (t)\in Y$ and let  $i$ be such that $z\in V_i.$ By continuity, let $c\in (a, b)$ be such that  $f\circ \alpha ([c,b))\subseteq V_i.$ Then $\alpha _{|}:(c,b)\to f^{-1}(V_i)\subseteq X$ is a definable curve in $f^{-1}(V_i)$ such that $f_{|}\circ \alpha _{|}:(c,b)\to V_i$ is completable. By hypothesis, $\alpha _{|}:(c,b)\to f^{-1}(V_i)\subseteq X$ is completable in $f^{-1}(V_i)$ and so $\lim_{t\to b^-}\alpha (t)$ exists in $X$ as required. \qed \\

From Corollary \ref{cor def proper4}  we obtain as in Corollary \ref{cor completion in B} the following analogue for definably proper.  In the case of o-minimal expansions of real closed fields this can be read off from  \cite[Chapter 10, (2.6) and (2.7)]{vdd}. \\

\begin{cor}\label{cor def completion}
 Let ${\mathbf B}$ be a  full subcategory of $\Df$. Suppose that the set of objects of ${\mathbf B}$ is:
 \begin{itemize}
 \item
 closed under taking  locally closed definable subspaces of objects of ${\mathbf B}$;
 \item
 is closed under taking cartesian products of objects of ${\mathbf B}$.
 \end{itemize}
 Then the following are equivalent:
\begin{enumerate}
\item
Every object $X$ of ${\mathbf B}$ is {\it definably completable in ${\mathbf B}$} i.e., there exists a definably compact space $X'$ in ${\mathbf B}$ together with a definable open immersion $i: X\hookrightarrow X'$ in ${\mathbf B}$ with  $i(X)$ dense in $X'$. Such $i: X\hookrightarrow X'$ is called a {\it definable completion} of $X$ in ${\mathbf B}$.

\item
Every  morphism $f:X\to Y$ in ${\mathbf B}$ is {\it definably completable in ${\mathbf B}$} i.e., there exists a commutative diagram
$$
\xymatrix{ X \ar [d]_f \ar[r]^-{i} & X' \ar[d]^{f'} \\
Y \ar[r]^-{j} & Y'}
$$
of  morphisms in ${\mathbf B}$ such that:  (i) $i:X\to X'$ is a definable  completion of $X$  in ${\mathbf B}$; (ii) $j$ is a definable  completion of $Y$ in ${\mathbf B}$.
\item
Every morphism $f:X\to Y$ in ${\mathbf B}$ has {\it a definable proper extension in ${\mathbf B}$} i.e., there exists a commutative diagram
$$
\xymatrix{ X \ar [rd]_f \ar[r]^\iota & P \ar[d]^{\overline{f}} \\
& Y}
$$
of morphisms in ${\mathbf B}$ such that $\iota$ is a definable  open immersion with $\iota (X)$ dense in $P$ and $\overline{f}$ is definably proper.
\end{enumerate}
If ${\mathbf B}=\Df$ we don't mention ${\mathbf B}$ and we talk of {\it definably completable, definable completion} and {\it definable proper extension}.\\
\end{cor}

\end{subsection}

\begin{subsection}{Definably proper and proper in $\Df$}\label{subsection def proper and proper Df}
Assuming  that ${\mathbb M}$ has definable Skolem functions,  below we will: (i)  show that a definably proper map between Hausdorff locally definably compact definable spaces  is the same a morphism proper in $\Df;$ (ii)  prove the definable analogue of the topological characterization of the notion of proper continuous maps (as closed maps with compact and Hausdorff fibers).\\ 

\begin{thm}\label{thm def proper2}
Suppose that ${\mathbb M}$ has definable Skolem functions. Let $X$ and $Y$ be  Hausdorff definable spaces with $Y$ locally definably compact. 
Let $f:X\to Y$ be a continuous definable map. Then the following are equivalent:
\begin{enumerate}
 \item
$f$ is  proper in $\Df$. 
\item
$f$ is definably proper.
\end{enumerate}
\end{thm}

\pf 
Assume (1). Let $\gamma : (a,b)\to X$ be a definable curve in $X$ and suppose that $f\circ \gamma :(a,b)\to Y$ is completable. By Theorem \ref{thm def proper and curves}, we need to show that $\gamma : (a,b)\to X$ is completable in $X$. By assumption $f\circ \gamma $ extends to a continuos definable map $g :[a,b]\to Y$. Consider a cartesian square of continuous definable maps
$$
\xymatrix{
(a,b) \ar@/_/[ddr] \ar@/^1.5pc/[drr]^{\gamma }
\ar@{-->}[dr]^-{\gamma  ' } \\
& X\times _Y [a,b] \ar[d]^-{f'} \ar[r]^-{g'}
& X \ar[d]^f \\
& [a,b] \ar[r]^g & Y }
$$
together with the continuous definable map $\gamma ':(a,b)\to X\times _Y[a,b]$ obtained from the maps $\gamma : (a,b)\to X$ and $(a,b) \hookrightarrow [a,b].$

Consider $\bar{\gamma '((a,b))}\subseteq X\times _Y[a,b].$
By assumption,  $f': X\times _Y[a,b]\to [a,b]$ is definably closed. So $f'(\bar{\gamma '((a,b))})$ is a closed definable subset of $[a,b]$. But $(a,b)=f'(\gamma '((a,b))\subseteq f'(\bar{\gamma '((a,b))})$ and so $f'(\bar{\gamma '((a,b))})=[a,b].$ Hence there are $u,v\in \bar{\gamma '((a,b))}$ such that $f'(u)=a$ and $f'(v)=b.$ Since $f'$ is the restriction of the projection $X\times [a,b]\to [a,b]$, $f'$ is definably open. Therefore, we have $\lim _{t\to a^+}\gamma '(t)=u$ and $\lim _{t\to b^-}\gamma '(t)=v$ respectively and $\gamma ': (a,b)\to X\times _Y[a,b]$ is completable in $X\times _Y[a,b].$ Thus $\gamma =g'\circ \gamma ' $ is completable in $X$ as required.

Assume  (2). Since  $f:X\to Y$ is separated in $\Df$ (Remark \ref{nrmk sep df2}), it is enough to show that $f$ is universally closed in $\Df$. For that it is enough to consider a cartesian square of continuous definable maps
$$
\xymatrix{X\times _YZ \ar[r]^{f'} \ar[d]^{g'} & Z \ar[d]^g \\
X \ar[r]^f &Y}
$$
and show that $f'$ is  definably closed (i.e. closed in $\Df$). 

Let $A\subseteq X\times _YZ$ be a closed definable subset and suppose that $f'(A)\subseteq Z$ is not closed. By the local almost everywhere curve selection (Theorem \ref{thm local aecs}) there is $z\in (\bar{f'(A)}\setminus f'(A))$ together with a definable curve  $\beta : (a,b)\to f'(A)\subseteq Y$  such that  $\lim _{t\to b^-}\beta  (t)=z.$   
By replacing $(a,b)$ by a smaller subinterval we may assume that $\lim _{t\to a^+}\beta (t)$ exists in $Z$, so $\beta $ is completable in $Z$. By definable Skolem functions, after replacing $(a,b)$ by a smaller subinterval, there exists a definable curve $\gamma :(a,b)\to X$ in $X$ such that for every $t\in (a,b)$ we have $\langle \gamma (t), \beta (t)\rangle \in A$.
Since $f\circ \gamma =g\circ \beta $ and $\beta $ is completable in $Z$, $g\circ \beta$ is completable in $Y$. Thus by (2) and Theorem \ref{thm def proper and curves}, $\gamma $ is completable in $X$ and  $\lim _{t\to b^-}\gamma (t)$ exists in $X$, call it $x$. If $\alpha =\langle \gamma , \beta \rangle :(a,b) \to X\times _YZ$, then $\lim _{t\to b^-}\alpha (t)=\langle x,z\rangle \in X\times _YZ$ and so $\langle x,z\rangle \in A$ because $A$ is closed. But then $z=f'(x,z)\in f'(A)$ which is absurd. \qed \\

Note that the assumption that $Y$ is locally definably compact is needed:\\

\begin{expl}\label{expl y lcompact}
{\em
Consider the setting of Example \ref{expl no aecs} and let $X=D,$ $Y=D\cup \{d\}$ and $f:X\to Y$ the inclusion. Then  $f:X\to Y$ is definably proper, $Y$ is not locally definably compact and $f$ is not definably closed.\\
}
\end{expl}



\begin{cor}\label{cor def normal comp3}
Suppose that ${\mathbb M}$ has definable Skolem functions. Let $X$ be a  Hausdorff  definable space. 
Then the following are equivalent:
\begin{enumerate}
\item
$X$ is definably compact.

\item
$X$ is complete in $\Df$. \\
\end{enumerate}
\end{cor}

The following is the definable analogue of the topological characterization of the notion of proper continuous maps (as closed maps with compact and Hausdorff fibers). A similar result appears in the semi-algebraic case (\cite[Theorem 12.5]{dk1}):\\

\begin{thm}\label{thm def proper fibers}
Suppose that ${\mathbb M}$ has definable Skolem functions. Let $X$ and $Y$ be  Hausdorff definable spaces with $Y$ locally definably compact. 
Let $f:X\to Y$ be a continuous definable map. Then the following are equivalent:
\begin{enumerate}
\item
 $f$ is definably proper.
\item
$f$ is definably closed  and has definably compact fibers.
\end{enumerate}
\end{thm}

\pf
Assume (1). Then $f:X\to Y$ has definably compact fibers and,  by Theorem \ref{thm def proper2}, $f$ is definably closed. 

Assume (2).  Let $K$ be a definably compact definable subset of $Y$. Let $\alpha :(a, b)\to f^{-1}(K)$ a definable curve in $f^{-1}(K)$. Suppose that $\lim _{t\to b^-}\alpha (t)$ does not exist in $f^{-1}(K)$. Then this limit does not exist in $X$ as well since $f^{-1}(K)$ is a closed definable subset of $X$ (by Corollary \ref{cor def comp2}, $K$ is closed). Therefore, if $d\in (a,b)$, then for every $e\in [d, b )$, $\alpha ([e, b))$ is a closed definable subset of $X$ contained in $f^{-1}(K)$. Indeed, we can first replace $a$ by $a'\in (a,b)$ if necessary so that $\alpha $ is injective and so $\alpha ((a,b))$ has a definable total order such that $\alpha $ is increasing; then if $\alpha ([e, b))$ is not closed, one can use the local almost everywhere curve selection (Theorem \ref{thm local aecs}) to obtain a definable curve $\delta : (a',b')\to \alpha ([e, b))$ with, say  $\lim _{t\to b'^-}\delta (t) \in {\rm cl}_X(\alpha ([e, b)))\setminus \alpha ([e, b));$ replacing $a'$ by some $a''\in  (a',b')$ if necessary, $\delta $ will be strictly increasing, but then we would have  $\lim _{t\to b^-}\alpha (t) =\lim _{t\to b'^-}\delta (t).$

 By assumption, for every $e\in [d, b ),$ $f\circ \alpha ([e, b))$ is then a closed definable subset of $Y$ contained in $K$. Since $K$ is definably compact, the limit $\lim _{t\to b^-}f\circ \alpha (t)$ exists in $K$, call it $c$. Hence, $c\in f\circ \alpha ([e, b))$ for every $e\in [d, b).$ Since the definable subset $\{t\in [d, b):f\circ \alpha (t)=c\}$ is a finite union of points and intervals, it follows that there is $d'\in [d,b)$ such that $f\circ \alpha (t)=c$ for all $t\in [d', b)$. Thus $\alpha ([d',b))\subseteq f^{-1}(c)\subseteq f^{-1}(K)$. Since $f^{-1}(c)$ is definably compact, the $\lim _{t\to b^-}\alpha (t)$ exists in $f^{-1}(K)$, which is absurd.
\qed \\

By Example \ref{expl y lcompact} the assumption that $Y$ is locally definably compact is needed.\\

\end{subsection}
\end{section}

\begin{section}{Invariance and comparison results}\label{section inv and comp}

\begin{subsection}{Definably proper in elementary extensions}\label{subsection def proper in ext}
Here ${\mathbb S}$ is an elementary extension of  ${\mathbb M}$ and we consider the  functor   
$$\Df \to \Df ({\mathbb S})$$
 from the category of definable spaces and continuous definable maps to the category of ${\mathbb S}$-definable spaces and continuous ${\mathbb S}$-definable maps.   This functor sends a definable space $X$ to the ${\mathbb S}$-definable space $X({\mathbb S})$ and sends a continuous definable map $f:X\to Y$ to the continuous ${\mathbb S}$-definable map $f^{{\mathbb S}}:X({\mathbb S})\to Y({\mathbb S})$. We show that: (i) $f$ is proper in $\Df$ if and only if $f^{{\mathbb S}}$ is proper in $\Df ({\mathbb S})$ (Theorem \ref{thm sep and proper in s}); (ii) if ${\mathbb M}$ has definable Skolem functions and $Y$ is  Hausdorff, then $f$ is  definably proper if and only if $f^{{\mathbb S}}$ is  ${\mathbb S}$-definably proper (Theorem \ref{thm def proper inv1}).\\

The following is easy and well known:\\

\begin{fact}\label{fact s def skolem1}
If ${\mathbb M}$ has definable Skolem functions, then ${\mathbb S}$ has definable Skolem functions. \\
\end{fact}

Since functor $\Df \to \Df ({\mathbb S})$  is a monomorphism from  the boolean algebra of definable subsets of  a definable space $X$ and  the boolean algebra of ${\mathbb S}$-definable  subsets of $X({\mathbb S})$ and it commutes with:
\begin{itemize}
\item
 the interior and closure operations;
 \item
 the  image and inverse image under (continuous) definable maps;
 \end{itemize}
we have:\\

\begin{lem} \label{lem def closed and in s}
Let $f:X\to Y$ a morphism in $\Df.$ Then the following are equivalent:
\begin{enumerate}
\item
$f$ is closed in $\Df $ (i.e. definably closed).
\item
$f^{{\mathbb S}}$ is closed in $\Df ({\mathbb S})$ (i.e. ${\mathbb S}$-definably closed).
\end{enumerate}
\end{lem}

\pf
Assume (1).  Let $A\subseteq X({\mathbb S})$ be a closed  ${\mathbb S}$-definable subset and suppose that $f^{{\mathbb S}}(A)$ is not a closed subset of $Y({\mathbb S})$. Then there is a uniformly definable family $\{A_t:t\in T\}$ of definable subsets of $X$ such that $A=A_s({\mathbb S})$ for some $s\in T({\mathbb S}).$ Since the property on $t$ saying that $A_t$ is closed is first-order, after replacing $T$ by a definable subset we may assume that for all $t\in T$, $A_t$ is a closed definable subset of $X$.  We also have that $\{f(A_t):t\in T\}$ is a uniformly definably family of definable subsets of $Y$ such that $f^{{\mathbb S}}(A)=f^{{\mathbb S}}(A_s({\mathbb S}))$.  Let $E$ be the  definable subset of $T$ of all $t$ such that $f(A_t)$ is not closed. Since $s\in E({\mathbb S})$, we have $E\neq \emptyset $ which is a contradiction since by assumption, for every $t\in T$, $f(A_t)$ is a closed definable subset of $Y$.

Assume (2).  Let $A\subseteq X$ be a closed definable subset. Then $A({\mathbb S})\subseteq X({\mathbb S})$ is a closed  ${\mathbb S}$-definable subset and by assumption,  $f(A)({\mathbb S})=f^{{\mathbb S}}(A({\mathbb S}))$ is a closed ${\mathbb S}$-definable subset of $Y({\mathbb S})$. So $f(A)$ is a closed definable subset of $Y$. \qed \\


Since functor $\Df \to \Df ({\mathbb S})$  sends open (resp. closed) definable immersion to open (resp. closed) $\mathbb{S}$-definable immersion and sends cartesian squares in $\Df$ to cartesian squares in $\Df ({\mathbb S})$ we have, using Lemma \ref{lem def closed and in s}: \\

\begin{thm} \label{thm sep and proper in s}
Let $f:X\to Y$ a morphism in $\Df.$ Then the following are equivalent:
\begin{enumerate}
\item
$f$ is proper (resp. separated) in $\Df .$
\item
$f^{{\mathbb S}}$ is proper (resp. separated) in $\Df ({\mathbb S}).$ \\
\end{enumerate}
\end{thm}

We also have:\\

\begin{thm}\label{thm def proper inv1}
Suppose that ${\mathbb M}$ has definable Skolem functions. Let $X$ and $Y$ be  definable spaces with $Y$ Hausdorff. 
Let $f:X\to Y$ be a continuous definable map. Then the following are equivalent:
\begin{enumerate}
\item
 $f$ is definably proper.
\item
$f^{{\mathbb S}}$ is ${\mathbb S}$-definably proper.
\end{enumerate}
\end{thm}

\pf
First note that  ${\mathbb S}$ has definable Skolem functions (Fact \ref{fact s def skolem1}) and  $Y({\mathbb S})$ is Hausdorff ${\mathbb S}$-definable spaces (since Hausdorff is a first-order property).  Using Corollary \ref{cor def comp2} and  Theorem \ref{thm def proper and curves}  in ${\mathbb M}$  and Corollary \ref{cor def comp2} and  Theorem \ref{thm def proper and curves} in ${\mathbb S},$ the result follows from the claim:

\begin{clm}\label{clm inv in s2}
The following are equivalent:
\begin{enumerate}
\item
There is a definable curve $\alpha :(a,b)\to  X$ such that $f\circ \alpha :(a,b)\to Y$ is completable in $Y$ but $\alpha $ is not completable in $X.$
\item
There is an ${\mathbb S}$-definable curve $\beta :(c,d)\to  X({\mathbb S})$ such that the  ${\mathbb S}$-definable curve  $f^{{\mathbb S}}\circ \beta :(c,d)\to Y({\mathbb S})$ is completable in $Y({\mathbb S})$ but $\beta $ is not completable in $X({\mathbb S}).$ \\
\end{enumerate}
\end{clm}

Assuming (1) then (2) holds with $(c,d)=(a,b)({\mathbb S})$ and $\beta =\alpha ^{{\mathbb S}}$ since ``$\alpha $ is continuous'',  ``$f\circ \alpha :(a,b)\to Y$ is completable in $Y$'' and  ``$\alpha $ is not completable in $X$'' are first-order properties.

Assume (2) then (1) also holds since ``$\beta $ is continuous'', ``$f^{{\mathbb S}}\circ \beta :(c,d)\to Y({\mathbb S})$ is completable in $Y({\mathbb S})$'' and ``$\beta $ is not completable in $X({\mathbb S})$'' are first-order properties in the parameters defining $\beta $ (together with $c$ and $d$).
\qed \\


The proof of Claim \ref{clm inv in s2} above actually shows:\\

\begin{cor}\label{cor def proper inv1}
Let $X$  be a 
  definable space. 
Then the following are equivalent:
\begin{enumerate}
\item
 $X$ is definably compact.
\item
$X({\mathbb S})$ is ${\mathbb S}$-definably compact.\\
\end{enumerate}
\end{cor}



\end{subsection}

\begin{subsection}{Definably proper in o-minimal expansions}\label{subsection proper in exp} 
Here  ${\mathbb S}$ is an o-minimal expansion of  ${\mathbb M}$ and we consider the  functor   
$$\Df \to \Df ({\mathbb S})$$
 from the category of definable spaces and continuous definable maps to the category of ${\mathbb S}$-definable spaces and continuous ${\mathbb S}$-definable maps.   This functor sends a definable space $X$ to the ${\mathbb S}$-definable space $X$ and sends a continuous definable map $f:X\to Y$ to the continuous ${\mathbb S}$-definable map $f:X\to Y$. We show that if ${\mathbb M}$ has definable Skolem functions, $X$ and $Y$ are  Hausdorff  and $Y$ is locally definably compact, then $f$ is definably proper if and only if $f^{{\mathbb S}}$ is  ${\mathbb S}$-definably proper and $f$ is proper in $\Df$ if and only if $f^{{\mathbb S}}$ is  proper in $\Df ({\mathbb S})$ (Theorem \ref{thm def proper inv3}).\\

\begin{fact}\label{fact s def skolem}
If ${\mathbb M}$ has definable Skolem functions, then ${\mathbb S}$ has definable Skolem functions.
\end{fact}

\pf
By Fact \ref{fact s def skolem1} we may assume that both ${\mathbb M}$ and ${\mathbb S}$ are $\omega $-saturated. In this case, by the (observations before the) proof of \cite[Chapter 6, (1,2)]{vdd} (see also Comment (1.3) there), ${\mathbb S}$ has definable Skolem functions if and only if  every nonempty ${\mathbb S}$-definable subset $X\subseteq M$ defined with parameters in $\bar{a}=a_1,\ldots , a_l$  has an element in ${\rm dcl}_{{\mathbb S}}(\bar{a}).$  

So let $X$ be an ${\mathbb S}$-definable subset of $M.$ By o-minimality, $X$ is a finite union of points $\{c_0, \ldots , c_m\}\subseteq M$ and open intervals $I_0, \ldots , I_n\subseteq M$ with end points  in $M\cup \{-\infty, +\infty \}$ with all the $c_i$'s and the endpoints of the $I_k$'s in ${\rm dcl}_{{\mathbb S}}(\bar{a}).$  So $X$ is definable over ${\rm dcl}_{{\mathbb S}}(\bar{a})$ using just equality and the order relation, hence $X$ is ${\mathbb M}$-definable. Since ${\mathbb M}$ has definable Skolem functions $X$ has a point in ${\rm dcl}_{{\mathbb M}}({\rm dcl}_{{\mathbb S}}(\bar{a})).$ Since ${\mathbb S}$ is an expansion of ${\mathbb M}$, we have ${\rm dcl}_{{\mathbb M}}({\rm dcl}_{{\mathbb S}}(\bar{a}))\subseteq {\rm dcl}_{{\mathbb S}}({\rm dcl}_{{\mathbb S}}(\bar{a}))={\rm dcl}_{{\mathbb S}}(\bar{a}).$

\qed \\

The shrinking lemma gives the following:\\

\begin{prop}\label{prop def comp and s comp}
Suppose that ${\mathbb M}$ has definable  Skolem functions. 
Let $X$ be a Hausdorff
definable space. Then the following are equivalent:
\begin{enumerate}
\item
$X$ is definably compact.
\item
 $X$ is   ${\mathbb S}$-definably compact.
 \end{enumerate}
\end{prop}

\pf
By Theorem \ref{thm def comp hausd normal}, $X$ is definably normal. Let  $(X_i, \phi _i)_{i\leq l}$ be the definable charts of $X.$  By the shrinking lemma, there are open definable subsets $V_i$ ($1\leq i \leq l$) and closed definable subsets $C_i$ ($1\leq i \leq l$) such that $V_i\subseteq C_i\subseteq X_i$ and $X=\cup \{C_i:i=1, \ldots , l\}.$  

Then  we have that $X$ is definably compact if and only if each $C_i$ is a definably compact definable subset of $X$ if and only if each  $\phi _i(C_i)$ is also a definably compact definable subset of $M^{n_i},$  and therefore, by \cite[Theorem 2.1]{ps}, if and only if each $\phi _i (C_i)$ is a closed and bounded definable subset of $M^{n_i}.$ 
Similarly we have that $X$ is ${\mathbb S}$-definably compact if and only if each $C_i$ is an ${\mathbb S}$-definably compact ${\mathbb S}$-definable subset of $X$ if and only if each  $\phi _i(C_i)$ is also an ${\mathbb S}$-definably compact ${\mathbb S}$-definable subset of $M^{n_i},$  and therefore, by \cite[Theorem 2.1]{ps} in ${\mathbb S}$, if and only if each $\phi _i (C_i)$ is a closed and bounded ${\mathbb S}$-definable subset of $M^{n_i}.$ 
Since ``closed'' and ``bounded'' are preserved  under going to ${\mathbb S}$ the result now follows.  \qed  \\

\begin{thm}\label{thm def proper inv3}
Suppose that ${\mathbb M}$ has definable Skolem functions. Let $X$ and $Y$ be  Hausdorff definable spaces with $Y$ locally definably compact.
Then the following are equivalent:
\begin{enumerate}
\item
$f$ is proper in $\Df$.
\item
 $f$ is definably proper.
\item
$f$ is ${\mathbb S}$-definably proper.
\item
$f$ is proper in $\Df ({\mathbb S}).$
\end{enumerate}
\end{thm}

\pf
First note that since $Y$ is locally definably compact, $Y({\mathbb S})$ is locally ${\mathbb S}$-definably compact.   By Theorem \ref{thm def proper2} in ${\mathbb M}$ and in ${\mathbb S}$ it is enough to show that  $f$ is definably proper if and only if $f$ is ${\mathbb S}$-definably proper. Using the fact that $Y$ is locally definably compact and Proposition \ref{prop def comp and s comp} one can show the later claim as in \cite[Chapter 6, (4.8) Exercise 2]{vdd} (see page 170 for the solution).
\qed \\




\end{subsection}

\begin{subsection}{Definably proper in topology}\label{subsection proper in top} 
Here  ${\mathbb M}$ is an o-minimal expansion of  the ordered set of real numbers and we consider the  functor   
$$\Df \to {\rm Top}$$
 from the category of definable spaces and continuous definable maps to the category of topological  spaces and continuous  maps.   We show that if ${\mathbb M}$ has definable Skolem functions, then for Hausdorff locally definably compact definable spaces definably proper is the same as  proper and proper in $\Df$ is the same as proper in ${\rm Top}$.\\

As before we have:

\begin{prop}\label{prop def comp and  comp}
Suppose that ${\mathbb M}$ has definable  Skolem functions. 
Let $X$ be a Hausdorff
definable space. Then the following are equivalent:
\begin{enumerate}
\item
$X$ is definably compact.
\item
 $X$ is   compact.
 \end{enumerate}
\end{prop}

\pf
Follow the proof of Proposition \ref{prop def comp and s comp} using the Heine-Borel theorem (a subset of ${\mathbb R}^n$ is compact if and only if it is closed and bounded) instead of \cite[Theorem 2.1]{ps}.
\qed \\

 A result  similar to the following  appears in the semi-algebraic case with a completely different proof (\cite[Theorem 9.11]{dk1}):\\

\begin{thm}\label{thm def proper inv4}
Suppose that ${\mathbb M}$ has definable Skolem functions. Let $X$ and $Y$ be  Hausdorff definable spaces with $Y$  locally definably compact.
Then the following are equivalent:
\begin{enumerate}
\item
$f$ is proper in $\Df$.
\item
 $f$ is definably proper.
\item
$f$ is proper.
\item
$f$ is proper in ${\rm Top}.$
\end{enumerate}
\end{thm}

\pf
First note that $Y$ is locally compact. Next recall that $f$ is proper if $f^{-1}(K)\subseteq X$ is a compact subset for every $K\subseteq Y$ compact subset and $f$ is proper in ${\rm Top}$  if it is separated and universally closed in the category ${\rm Top}$ of topological spaces. Also, it is well know that $f$ is proper if and only if $f$ is proper in ${\rm Top}$ if and only if $f$ is closed and has compact fibers (see \cite[Chapter 1, $\S$10, Theorem 1]{bki}). 

By Theorem \ref{thm def proper2}  it is enough to show that  $f$ is definably proper if and only if $f$ is  proper. Using Theorem \ref{thm def proper fibers} and Proposition \ref{prop def comp and comp} one can show the later claim as in \cite[Chapter 6, (4.8) Exercise 3]{vdd} (see page 170 for the solution).
\qed \\

\end{subsection}

\end{section}

\begin{section}{Definably compact, definably proper and definable types}\label{section def comp, def proper and def types}
Here  we show that definable compactness  of Hausdorff definable spaces in o-minimal structures with definable Skolem functions can also be characterized by existence of limits of definable types - extending a similar result in the affine case (\cite[Remark 4.2.15]{HrLo}). The corresponding characterization for definably proper maps between Hausdorff, locally definably compact definable spaces is also given (Theorem \ref{thm def proper1}).\\

Let $X$ be a definable space. A {\it type on $X$} is an ultrafilter $\alpha $ of definable subsets of $X$. A type $\alpha $ on $X$ is a {\it definable type on $X$} if for every uniformly definable family $\{F_t\}_{t\in T}$ of definable subsets of $X$, with $T\subseteq M^n$ for some $n$, there is a definable subset $T(\alpha )\subseteq T$ such that $F_t\in \alpha $ if and only if $t\in T(\alpha ).$ 

If $\alpha $ is a type on $X$ and $x\in X$, we say that {\it $x$ is a limit of $\alpha $},  if for every open definable subset $U$ of $X$ such that $x\in U$ we have $U\in \alpha .$\\

For affine definable spaces existence of limits of definable types gives another criteria for definable compactness (see  \cite[Remark 4.2.15]{HrLo}). Since the proof is not written down in \cite{HrLo}, for convenience, we include the details.\\

\begin{fact}\label{fact def comp types}
Let $Z\subseteq M^n$ be a definable set. Then the following are equivalent:
\begin{enumerate}
\item
$Z$ is closed and bounded (i.e., definably compact). 
\item
Every definable type on $Z$ has a limit in $Z.$
\end{enumerate}
\end{fact}

\pf
Assume (1). Let $\alpha $ be a definable type on $Z.$ For each $i=1, \ldots, n,$ let $\pi _i:M^n\to M$ be the projection onto the $i$-coordinate and let $Z_i=\pi _i(Z)$ and $\alpha _i=\tilde{\pi _i}(\alpha )$ (the definable type on $Z_i$ determined by the collection of definable subsets $\{A\subseteq Z_i: \pi _i^{-1}(A)\in \alpha \}$). By \cite[Lemma 2.3]{MarStein}, each $\alpha _i$ is not a cut  and so, since each $Z_i$ is bounded,  there is $a_i\in M$ such that $\alpha _i$ is either determined by $x=a_i$ or  $\{b<x<a_i: b\in M, \,\,b<a_i\}$ or $\{a_i<x<b: b\in M, \,\,a_i<b\}.$ In either case $a_i$ is the limit of $\alpha _i$ in $M.$  Clearly $a=\langle a_1,\ldots , a_n \rangle\in M^n$ is the limit of $\alpha $ in $M^n$ and, since $Z$ is closed, $a\in Z$ as required.

Assume (2). Let $\alpha :(a,b)\to Z$ be a definable curve. Then the collection of definable subsets $\{\alpha ([t,b)): t\in (a,b)\}$ of $Z$ determines a  type $\beta $ on $Z$ such that $S\in \beta $ if and only if $\alpha ([t,b))\subseteq S$ for some $t\in (a,b).$ This type $\beta $ is definable since for any uniformly definable family $\{F_l\}_{l\in L}$ of definable subsets of $Z$ we have $\{l\in L: F_l\in \beta \}=\{l\in L: \exists t\in (a,b)(\alpha ([t,b))\subseteq F_l)\}.$
By hypothesis, $\beta $ has a limit $z$ in $Z$. This $z\in Z$ is also  the limit $\lim _{t\to b^-}\alpha (t)$ since for every  $d\in D(z)$ we have $U(z,d)\in \beta ,$ so $\alpha ([t,b))\subseteq U(z,d)$ for some $t\in (a,b).$
\qed \\

We can use the shrinking lemma to extend this result to non affine Hausdorff definable spaces:\\

\begin{thm}\label{thm def comp norm types}
Suppose that ${\mathbb M}$ has definable Skolem functions. Let $X$ be a Hausdorff definable space. Then the following are equivalent:
\begin{enumerate}
\item
$X$ is definably compact.
\item
Every definable type on $X$ has a limit in $X$.
\end{enumerate}
\end{thm}

\pf
By Theorem \ref{thm def comp hausd normal}, $X$ is  definably normal. Let  $(X_i, \theta _i)_{i\leq k}$ be the definable charts of $X$ with $\theta _i(X_i)\subseteq M^{n_i}.$ By the shrinking lemma  there are definable open subsets $V_i$ and definable closed subsets $C_i$ of $X$ ($1\leq i\leq n$) with $V_i\subseteq C_i\subseteq X_i$ and $X=\cup \{C_i:i=1,\dots, n\}$. We have that: (i) $X$ is definably compact if and only if each $C_i$ is definably compact; (ii) every definable type on $X$ has limit in $X$ if and only if for each $i$, every definable type on $C_i$ has a limit in $C_i.$ Since $\theta _{i|}:C_i\to \theta _i(C_i)\subseteq M^{n_i}$ is a definable homeomorphism the result now follows  by Fact \ref{fact def comp types}.
\qed \\

We also have the following definable types criterion for definably proper:\\

\begin{thm}\label{thm def proper1}
Suppose that ${\mathbb M}$ has definable Skolem functions. Let $X$ and $Y$ be  Hausdorff definable spaces with $Y$ locally definably compact. 
Let $f:X\to Y$ be a continuous definable map. Then the following are equivalent:
\begin{enumerate}
\item
 $f$ is definably proper.
\item
For every definable type $\alpha $ on $X$, if $\tilde{f}(\alpha )$ has a limit in $Y$, then $\alpha $ has a limit in $X.$
\end{enumerate}
\end{thm}

\pf
Assume (1). Let $\alpha $ be a definable type on $X$ such that $\tilde{f}(\alpha )$ has a limit in $Y$, say $\lim \tilde{f}( \alpha ) =y\in Y$.  Since $Y$ is locally definable compact, there is a definable open neighborhood $V$ of $y$ in $Y$ such that $\bar{V}$ is definably compact. So, $f^{-1}(\bar{V})$ is a definably compact definable subset of $X$  and $\alpha $ is a definable type on $f^{-1}(\bar{V})$. But then by  Theorem \ref{thm def comp norm types}  $\alpha  $ has a limit  in $f^{-1}(\bar{V})$, hence in $X$.

Assume (2). Suppose that $f$ is not definably proper. Then there is a definably compact definable subset $K$ of $Y$ such that $f^{-1}(K)$ is not a definably compact definable subset of $X$. Thus by  Theorem \ref{thm def comp norm types}  there is a definable type $\alpha  $ on $f^{-1}(K)$ which does  not have a limit  in $f^{-1}(K)$. Since $f^{-1}(K)$ is closed (by Corollary \ref{cor def comp2}, $K$ is closed), $\alpha  $ does  not have a limit in $X$. But $\tilde{f}(\alpha )$ is a definable type on $K\subseteq Y$ and has a limit  by Theorem \ref{thm def comp norm types}, which contradicts (2). 
 \qed \\

The following was observed in \cite[Remark 4.2.15]{HrLo} in the affine case but the same proof works. \\

\begin{fact}\label{fact nonclosed types}
Suppose that ${\mathbb M}$ has definable Skolem functions.  Let $X$ be a definable space and $C\subseteq X$ a definable subset which is not closed. If $x\in \bar{C}\setminus C$ then there is  a definable type $\alpha $ on $C$ such that $x$ is a limit of  $\alpha .$
\end{fact}

\pf
Consider  the definable set $D(x)$ with the relation $\preceq $ (a definable downwards directed order). By \cite[Lemma 4.2.18]{HrLo} (or \cite[Lemma 2.19]{Hr04}) there is a definable type $\beta $ on $D(x)$ such that for every $d\in D(x)$ we have $\{d' \in D(x): d'\preceq d\}\in \beta .$ 

Since $x\in \bar{C},$ for every $d\in D(x)$ we have that $U(x,d)\cap C\neq \emptyset .$ By definable Skolem functions, there is a definable map $h:D(x)\to C$ such that for every $d\in D(x)$ we have $h(d)\in U(x,d)\cap C.$ Let $\alpha =\tilde{h}(\beta )$ be the definable type on $C$ determined by the collection of definable subsets $\{A\subseteq C: h^{-1}(A)\in \beta \}.$ Clearly,  $x$ is a  limit of $\alpha. $
\qed \\

By Example \ref{expl y lcompact} and Fact \ref{fact nonclosed types}, in Theorem \ref{thm def proper1} the assumption that $Y$ is locally definably compact is needed. Note that, by the same example,  this observation applies also if one replaces the Peterzil-Steinhorn definition of definable compact (using definable curves - \cite{ps}) by the Hrushovski-Loeser definition of definable compact (using definable types - \cite{HrLo}). \\

\end{section}

\end{document}